\documentclass[12pt]{amsart}
\usepackage{amsmath,amssymb,amscd,latexsym,eufrak}

\newcommand{\svskip}{\vspace{3mm}}
\newcommand{\C}{{\mathbb C}}
\newcommand{\Q}{{\mathbb Q}}
\newcommand{\Z}{{\mathbb Z}}
\newcommand{\BP}{{\mathbb P}}

\newcommand{\gm}{\EuFrak{m}}

\newcommand{\QED}{{\unskip\nobreak\hfil\penalty50\quad\null\nobreak\hfil
{$\Box$}\parfillskip0pt\finalhyphendemerits0\par\medskip}}
\newcommand{\Proof}{\noindent{\bf Proof.}\quad}
\newcommand{\A}{{\mathbb A}}
\renewcommand{\H}{{\rm H}}
\newcommand{\N}{{\rm N}}
\newcommand{\Spec}{{\rm Spec}\:}

\newcommand{\lto}{\longrightarrow}

\newcommand{\Pic}{{\rm Pic}\:}

\newcommand{\Ker}{{\rm Ker}\:}
\newcommand{\Der}{{\rm Der}}

\newcommand{\Hom}{{\rm Hom}}
\newcommand{\Ext}{{\rm Ext}}
\newcommand{\SHom}{{\mathcal H}{\it om}}

\newcommand{\Hilb}{{\rm Hilb}}
\newcommand{\Alb}{{\rm Alb}}
\newcommand{\id}{{\rm id}}

\newcommand{\coim}{{\rm Coim}\:}

\newcommand{\lkd}{\ol{\kappa}}
\newcommand{\quot}{/\!/}
\newcommand{\SD}{{\mathcal D}}
\newcommand{\SF}{{\mathcal F}}
\newcommand{\SH}{{\mathcal H}}
\newcommand{\SO}{{\mathcal O}}

\newcommand{\ST}{{\mathcal T}}
\newcommand{\SW}{{\mathcal W}}
\newcommand{\SX}{{\mathcal X}}
\newcommand{\wt}{\widetilde}
\newcommand{\wh}{\widehat}
\newcommand{\ol}{\overline}

\newcommand{\ML}{{\rm ML}}

\newcommand{\dps}[1]{\displaystyle{#1}}

\newtheorem{thm}{Theorem}[section]
\newtheorem{lem}[thm]{Lemma}

\newtheorem{cor}[thm]{Corollary}
\newtheorem{remark}[thm]{Remark}
\newtheorem{example}[thm]{Example}
\newtheorem{prob}[thm]{Problem}

\newtheorem{claim}{\sc Claim}
\newtheorem{conj}[thm]{Conjecture}
\newcommand{\st}[1]{\stackrel{{#1}}{\longrightarrow}}

\begin{document}
\title[Deformations of $\A^1$-fibrations]
{Deformations of $\A^1$-fibrations}
\author{R.V. Gurjar, K. Masuda and M. Miyanishi}
\date{August 28, 2013}

\thanks{The second and third authors are supported by Grant-in-Aid for Scientific Research (C), No. 22540059 
and (B), No. 24340006, JSPS} 
\address{School of Mathematics, Tata Institute of Fundamental Research \\
Homi Bhabha Road, Mumbai 400 001, India}
\email{gurjar@math.tifr.res.in}
\svskip

\address{School of Science \& Technology, Kwansei Gakuin University \\
2-1 Gakuen, Sanda 669-1337, Japan}
\email{kayo@kwansei.ac.jp}
\svskip

\address{Research Center for Mathematical Sciences \\
Kwansei Gakuin University\\
2-1 Gakuen, Sanda 669-1337, Japan}
\email{miyanisi@kwansei.ac.jp}

\maketitle

\begin{abstract}
Let $B$ be an integral domain which is finitely generated over a subdomain $R$ and let $D$ be an $R$-derivation on $B$ 
such that the induced derivation $D_\gm$ on $B\otimes_RR/\gm$ is locally nilpotent for every maximal ideal $\gm$. We ask 
if $D$ is locally nilpotent. Theorem \ref{Theorem 1.1} asserts that this is the case if $B$ and $R$ are affine domains. 
We next generalize the case of $G_a$-action treated in Theorem \ref{Theorem 1.1} to the case of $\A^1$-fibrations and consider 
the log deformations of affine surfaces with $\A^1$-fibrations. The case of $\A^1$-fibrations of affine type behaves nicely 
under log deformations, while the case of $\A^1$-fibrations of complete type is more involved (see Dubouloz-Kishimoto \cite{DK}). 
As a corollary, we prove the generic triviality of $\A^2$-fibration over a curve and generalize this result to the case of 
affine pseudo-planes of $\ML_0$-type under a suitable monodromy condition.
\end{abstract}

\section*{Introduction}

An $\A^1$-fibration $\rho : X \to B$ on a smooth affine surface $X$ to a smooth curve $B$ is given as the quotient 
morphism of a $G_a$-action if the parameter curve $B$ is an affine curve (see \cite{GKM}). Meanwhile, it is not so if 
$B$ is a complete curve. When we deform the surface $X$ under a suitable setting (log deformation), our question is 
if the neighboring surfaces still have $\A^1$-fibrations of affine type or of complete type according to the type of 
the $\A^1$-fibration on $X$ being affine or complete. Assuming that the neighboring surfaces have $\A^1$-fibrations, 
the propagation of the type of $\A^1$-fibration is proved in Lemma \ref{Lemma 2.2}, whose proof reflects the structure 
of the boundary divisor at infinity of an affine surface with $\A^1$-fibration. The stability of the boundary divisor under 
small deformations, e.g., the stability of the weighted dual graphs has been discussed in topological methods (e.g., \cite{N}). 
Furthermore, if such property is inherited by the neighboring surfaces, we still ask if the ambient threefold has an 
$\A^1$-fibration or equivalently if the generic fiber has an $\A^1$-fibration.  

The answer to this question is subtle. We consider first in the section one the case where each of the fiber surfaces 
of the deformation has an $\A^1$-fibration of affine type induced by a global vector field on the ambient threefold. 
This global vector field is in fact given by a locally nilpotent derivation (Theorem \ref{Theorem 1.1}). If the 
$\A^1$-fibrations on the fiber surfaces are of affine type, we can show (Theorem \ref{Theorem 2.8}), with the absence of 
monodromies of boundary components, that there exists an $\A^1$-fibration on the ambient threefold such that the 
$\A^1$-fibration on each general fiber surface is induced by the global one up to an automorphism of the fiber surface. 
The proof of Theorem \ref{Theorem 2.8} depends on Lemma \ref{Lemma 2.2} which we prove by observing the behavior of 
the boundary rational curves. This is done by the use of Hilbert scheme (see \cite{Kollar})and by killing monodromies by 
\'etale finite changes of the base curve. 

As a consequence, we can prove the generic triviality of an $\A^2$-fibration over a curve. Namely, if $f : Y \to T$ is a 
smooth morphism from a smooth affine threefold to a smooth affine curve such that the fiber over every {\em closed} point 
of $T$ is isomorphic to the affine plane $\A^2$, then the generic fiber of $f$ is isomorphic to $\A^2$ over the function 
field $k(T)$ of $T$ and $f$ is an $\A^2$-bundle over an open set of $T$ (see Theorem \ref{Theorem 2.10}). This fact, 
together with a theorem of Sathaye \cite{S}, shows that $f$ is an $\A^2$-bundle over $T$ in the Zariski topology.

The question on the generic triviality is also related to a question on the triviality of a $k$-form of a surface with 
an $\A^1$-fibration (see Problem \ref{Problem 2.13}). In the case of an $\A^1$-fibration of complete type, 
the answer is negative by Dubouloz-Kishimoto \cite{DK} (see Theorem \ref{Theorem 5.1}).

Theorem \ref{Theorem 2.10} was proved by our predecesors Kaliman-Zaidenberg \cite{KZ} in a more comprehensive way and 
without assuming that the base is a curve. The idea in our first proof of Theorem \ref{Theorem 2.10} is of more algebraic 
nature and consists of using the existence of a locally nilpotent derivation on the coordinate ring of $Y$ and 
the second proof of using the Ramanujam-Morrow graph of the normal minimal completion of $\A^2$ was already used in 
\cite{KZ}. The related results are also discussed in the article \cite{Russell}. We cannot still avoid the use of a 
theorem of Kambayashi \cite{Kamb} on the absence of separable forms of the affine plane.

Some of the algebro-geometric arguments using Hilbert scheme in the section two can be replaced by topological arguments 
using Ehresmann's theorem which might be more appreciated than the use of the Hilbert scheme. But they are restricted to the 
case of small deformations. This is done in the section three. 

In the section four, we extend the above result on the generic triviality of an $\A^2$-fibration over a curve by 
replacing $\A^2$ by an affine pseudo-plane of $\ML_0$-type which has properties similar to $\A^2$, e.g., the boundary 
divisor for a minimal normal completion is a linear chain of rational curves. But we still need a condition on the 
monodromy. An affine pseudo-plane, not necessarily of $\ML_0$-type, is a $\Q$-homology plane, and we note that 
Flenner-Zaidenberg \cite{FZ} made a fairly exhaustive consideration for the log deformations of $\Q$-homology planes.

In the final section five, we observe the case of $\A^1$-fibration of complete type and show by an example of 
Dubouloz-Kishimoto \cite{DK} that the ambient threefold does not have an $\A^1$-fibration. But it is still plausible 
that the ambient threefold is affine-uniruled in the stronger sense that the fiber product of the ambient deformation 
space by a suitable lifting of the base curve has a global $\A^1$-fibration. But this still remains open. 

We use two notations for the intersection of (not necessarily irreducible) subvarieties $A, B$ of codimension one 
in an ambient threefold. Namely, $A\cap B$ is the intersection of two subvarieties, and $A\cdot B$ is the 
intersection of effective divisors. In most cases, both are synonymous. 

As a final remark, we note that a preprint of Flenner-Kaliman-Zaidenberg \cite{FKZ} recently uploaded on the web treats 
also deformations of surfaces with $\A^1$-fibrations. 

The referees pointed out several flaws in consideration of monodromies in the preliminary versions of the article. In 
particular, we are indebted for Example \ref{Example 2.6} to one of the referees. The authors are very grateful to 
the referees for their comments and advice.

\section{Triviality of deformations of locally nilpotent derivations}

Let $k$ be an algebraically closed field of characteristic zero which we fix as the ground field. Let 
$Y=\Spec B$ be an irreducible affine algebraic variety. We define the tangent sheaf $\ST_{Y/k}$ as 
$\SHom_{\SO_Y}(\Omega^1_{Y/k}, \SO_Y)$. A regular vector field on $Y$ is an element of $\Gamma(Y,\ST_{Y/k})$. 
A regular vector field $\Theta$ on $Y$ is identified with a derivation $D$ on $B$ via isomorphisms
\[
\Gamma(Y,\ST_{Y/k})\cong \Hom_B(\Omega^1_{B/k},B)\cong \Der_k(B,B).
\]
We say that $\Theta$ is {\em locally nilpotent} if so is $D$. In the first place, we are interested in finding a necessary 
and sufficient condition for $D$ to be locally nilpotent. Suppose that $Y$ has a fibration $f : Y \to T$. A natural question 
is to ask whether $D$ is locally nilpotent if the restriction of $D$ on each closed fiber of $f$ is locally nilpotent.
The following result shows that this is the case\footnote{The result is also remarked in \cite[Remark 13]{DK}.}.

\begin{thm}\label{Theorem 1.1}
Let $Y=\Spec B$ and $T=\Spec R$ be irreducible affine varieties defined over $k$ and let $f : Y \to T$ be a dominant 
morphism such that general fibers are irreducible and reduced. We consider $R$ to be a subalgebra of $B$. Let $D$ be an 
$R$-trivial derivation of $B$ such that, for each closed point $t \in T$, the restriction $D_t=D\otimes_RR/\gm$ is a
locally nilpotent derivation of $B\otimes_RR/\gm$, where $\gm$ is the maximal ideal of $R$ corresponding to $t$. 
Then $D$ is locally nilpotent. 
\end{thm}

We need some preliminary results. We retain the notations and assumptions in the above theorem.

\begin{lem}\label{Lemma 1.2}
There exist a finitely generated field extension $k_0$ of the prime field $\Q$ which is a subfield of the ground 
field $k$, geometrically integral affine varieties $Y_0=\Spec B_0$ and $T_0=\Spec R_0$, a dominant morphism 
$f_0 : Y_0 \to T_0$ and an $R_0$-trivial derivation $D_0$ of $B_0$ such that the following conditions are 
satisfied:
\begin{enumerate}
\item[(1)]
$Y_0, T_0, f_0$ and $D_0$ are defined over $k_0$.
\item[(2)]
$Y=Y_0\otimes_{k_0}k, T=T_0\otimes_{k_0}k, f=f_0\otimes_{k_0}k$ and $D=D_0\otimes_{k_0}k$.
\item[(3)]
$D_0$ is locally nilpotent if and only if so is $D$.
\end{enumerate}
\end{lem}

\Proof
Since $B$ and $R$ are integral domains finitely generated over $k$, write $B$ and $R$ as the residue rings of certain 
polynomial rings over $k$ modulo the finitely generated ideals. Write $B=k[x_1, \ldots,x_r]/I$ and $R=k[t_1, \ldots,t_s]/J$. 
Furthermore, the morphism $f$ is determined by the images $f^*(\eta_j)=\varphi_j(\xi_1, \ldots, \xi_r)$ in $B$, where 
$\xi_i=x_i \pmod{I}$ and $\eta_j= t_j \pmod{J}$. Adjoin to $\Q$ all coefficients of the finite generators of $I$ and $J$ 
as well as the coefficients of the $\varphi_j$ to obtain a subfield $k_0$ of $k$. Let $B_0=k_0[x_1,\ldots, x_r]/I_0$ 
and $R_0=k_0[t_1,\ldots,t_s]/J_0$, where $I_0$ and $J_0$ are respectively the ideals in $k_0[x_1,\ldots,x_r]$ and 
$k_0[t_1,\ldots,t_s]$ generated by the same generators of $I$ and $J$. Furthermore, define the homomorphism $f_0^*$ by the 
assignment $f_0^*(\eta_j)=\varphi_j(\xi_1, \ldots,\xi_r)$. Let $Y_0=\Spec B_0$, $T_0=\Spec R_0$ and let $f_0 : Y_0 \to T_0$ 
be the morphism defined by $f_0^*$. The derivation $D$ corresponds to a $B$-module homomorphism $\delta : \Omega^1_{B/R} \to B$.
Since $\Omega^1_{B/R}=\Omega^1_{B_0/R_0}\otimes_{k_0}k$, we can enlarge $k_0$ so that there exists a $B_0$-homomorphism 
$\delta_0 : \Omega^1_{B_0/R_0} \to B_0$ satisfying $\delta=\delta_0\otimes_{k_0}k$. Let $D_0=\delta_0\cdot d_0$, where 
$d_0 : B_0 \to \Omega^1_{B_0/R_0}$ is the standard differentiation. Then we have $D=D_0\otimes_{k_0}k$. 

Let $\Phi_0 : B_0 \to B_0[[u]]$ be the $R_0$-homomorphism into the formal power series ring in $t$ over $B_0$ defined by
\[
\Phi_0(b_0)=\sum_{i \geq 0}\frac{1}{i!}D_0^i(b_0)u^i\ .
\]
Let $\Phi : B \to B[[u]]$ be the $R$-homomorphism defined in a similar fashion. Then $\Phi_0$ and $\Phi$ are determined 
by the images of the generators of $B_0$ and $B$. Since the generators of $B_0$ and $B$ are the same, we have $\Phi=
\Phi_0\otimes_{k_0}k$. The derivation $D_0$ is locally nilpotent if and only if $\Phi_0$ splits via the polynomial 
subring $B_0[u]$ of $B_0[[u]]$. This is the case for $D$ as well. Since $\Phi_0$ splits via $B_0[u]$ if and only if 
$\Phi$ splits via $B[u]$, $D_0$ is locally nilpotent if and only if so is $D$.
\QED 

\begin{lem}\label{Lemma 1.3}
Let $k_1$ be the algebraic closure of $k_0$ in $k$. Let $Y_1=\Spec B_1$ with $B_1=B_0\otimes_{k_0}k_1$, $T_1= 
\Spec R_1$ with $R_1=R_0\otimes_{k_0}k_1$ and $f_1=f_0\otimes_{k_0}k_1$. Let $D_1=D_0\otimes_{k_0}k_1$. Then the 
following assertions hold.
\begin{enumerate}
\item[(1)]
Let $t_1$ be a closed point of $T_1$. Then the restriction of $D_1$ on the fiber $f_1^{-1}(t_1)$ is locally nilpotent.
\item[(2)]
$D_1$ is locally nilpotent if and only if so is $D$.
\end{enumerate}
\end{lem}

\Proof
(1)\ Let $t$ be the unique closed point of $T$ lying over $t_1$ by the projection morphism $T \to T_1$, where 
$R=R_1\otimes_{k_1}k$. (If $\gm_1$ is the maximal ideal of $R_1$ corresponding to $t_1$, $\gm_1\otimes_{k_1}k$ is 
the maximal ideal of $R$ corresponding to $t$.) Then $F_t=f^{-1}(t)=f_1^{-1}(t_1)\otimes_{k_1}k$ and the restriction 
$D_t$ of $D$ onto $F_t$ is given as $D_{1,t_1}\otimes_{k_1}k$, where $D_{1,t_1}$ is the restriction of $D_1$ onto 
$f_1^{-1}(t_1)$. We consider also the $R$-homomorphism $\Phi : B \to B[[u]]$ and the $R_1$-homomorphism $\Phi_1 : 
B_1 \to B_1[[u]]$. As above, let $\gm$ and $\gm_1$ be the maximal ideals of $R$ and $R_1$ corresponding to $t$ and $t_1$. 
Then $D_t$ gives rise to the $R/\gm$-homomorphism $\Phi\otimes_RR/\gm : B\otimes_RR/\gm \to (B\otimes_RR/\gm)[[u]]$. 
Similarly, $D_{1,t_1}$ gives rise to the $R_1/\gm_1$-homomorphism $\Phi_1\otimes_{R_1}R_1/\gm_1 : 
B_1\otimes_{R_1}R_1/\gm_1 \to (B_1\otimes_{R_1}R_1/\gm_1)[[u]]$, where $R/\gm=k$ and $R_1/\gm_1=k_1$. Then 
$\Phi\otimes_RR/\gm=(\Phi_1\otimes_{R_1}R_1/\gm_1)\otimes_{k_1}k$. Hence $\Phi\otimes_RR/\gm$ splits via 
$(B\otimes_RR/\gm)[u]$ if and only if $\Phi_1\otimes_{R_1}R_1/\gm_1$ splits via $(B_1\otimes_{R_1}R_1/\gm_1)[u]$. 
Hence $D_{1,t}$ is locally nilpotent as so is $D_t$.

(2)\ The same argument as above using the homomorphism $\Phi$ can be applied.
\QED

The field $k_0$ can be embedded into the complex field $\C$ because it is a finitely generated field extension of $\Q$.
Hence we can extend the embedding $k_0 \hookrightarrow \C$ to the algebraic closure $k_1$. Thus $k_1$ is viewed as a 
subfield of $\C$. Then Lemma \ref{Lemma 1.3} holds if one replaces the extension $k/k_1$ by the extension $\C/k_1$. 
Hence it suffices to prove Theorem \ref{Theorem 1.1} with an additional hypothesis $k=\C$.

\begin{lem}\label{Lemma 1.4}
Theorem \ref{Theorem 1.1} holds if $k$ is the complex field $\C$.
\end{lem}

\Proof
Let $Y(\C)$ be the set of closed points which we view as a complex analytic space embedded into a complex affine space $\C^N$
as a closed set. Consider the Euclidean metric on $\C^N$ and the induced metric topology on $Y(\C)$. Then $Y(\C)$ is a 
complete metric space. 

Let $b$ be a nonzero element of $B$. For a positive integer $m$, define a Zariski closed subset $Y_m(b)$ of $Y(\C)$ by
\[
Y_m(b)=\{Q \in Y(\C) \mid D^m(b)(Q)=0\}\ .
\]
Since $Q$ lies over a closed point $t$ of $T(\C)$ and $D_t$ is locally nilpotent on $f^{-1}(t)$ by the hypothesis, we have 
\[
f^{-1}(t) \subset \bigcup_{m>0}Y_m(b)\ .
\]
This implies that $Y(\C)=\dps{\bigcup_{m>0}Y_m(b)}$. We claim that $Y(\C)=Y_m(b)$ for some $m > 0$.  In fact, this follows 
by Baire category theorem, which states that if the $Y_m(b)$ are all proper closed subsets, its countable union cannot 
cover the uncountable set $Y(\C)$. If $Y(\C)=Y_m(b)$ for some $m > 0$ then $D^m(b)=0$. This implies that $D$ is locally 
nilpotent on $B$.

One can avoid the use of Baire category theorem in the following way. Suppose that $Y_m(b)$ is a proper closed subset 
for every $m > 0$. Let $H$ be a general hyperplane in $\C^N$ such that the section $Y(\C)\cap H$ is irreducible, 
$\dim Y(\C)\cap H=\dim Y(\C)-1$, and $Y(\C)\cap H=\dps{\bigcup_{m>0}(Y_m(b)\cap H)}$ with $Y_m(b)\cap H$ a proper 
closed subset of $Y(\C)\cap H$ for every $m > 0$. We can further take hyperplane sections and find a general linear 
subspace $L$ in $\C^N$ such that $Y(\C)\cap L$ is an irreducible curve and $Y(\C)\cap L=\dps{\bigcup_{m>0}(Y_m(b)\cap L)}$, 
where $Y_m(b)\cap L$ is a proper Zariski closed subset. Hence $Y_m(b)\cap L$ is a finite set, and $\dps{\bigcup_{m>0}(Y_m(b)\cap L)}$
is a countable set, while $Y(\C)\cap L$ is not a countable set. This is a contradiction. Thus $Y(\C)=Y_m(b)$ for some 
$m > 0$.
\QED

Let $D$ be a $k$-derivation on a $k$-algebra $B$. It is called {\em surjective} if $D$ is so as a $k$-linear mapping. 
The follwoing result is a consequence of Theorem \ref{Theorem 1.1}

\begin{cor}\label{Corollary 1.5}
Let $Y=\Spec B$, $T=\Spec R$ and $f : Y \to T$ be the same as in Theorem \ref{Theorem 1.1}. Let $D$ be an $R$-derivation 
of $B$ such that $D_t$ is a surjective $k$-derivation for every closed point $t \in T$. Assume further that the relative 
dimension of $f$ is one. Then $D$ is a locally nilpotent derivation and $f$ is an $\A^1$-fibration.
\end{cor}

\Proof
Let $t$ be a closed point of $T$ such that the fiber $f^{-1}(t)$ is irreducible and reduced. 
By \cite[Theorem 1.2 and Proposition 1.7]{GKM2}, the coordinate ring $B\otimes R/\gm$ of $f^{-1}(t)$ is a polynomial ring 
$k[x]$ in one variable and $D_t=\partial /\partial x$, where $\gm$ is the maximal ideal of $R$ corresponding to $t$. 
Then $D_t$ is locally nilpotent. Taking the base change of $f : Y \to T$ by $U \hookrightarrow T$ if necessary, where 
$U$ is a small open set of $T$, we may assume that $D_t$ is locally nilpotent for every closed point $t$ of $T$. By 
Theorem \ref{Theorem 1.1}, the derivation $D$ is locally nilpotent and hence $f$ is an $\A^1$-fibration.
\QED

\section{Deformations of $\A^1$-fibrations of affine type}

In the present section, we assume that the ground field $k$ is the complex field $\C$. Let $X$ be an affine algebraic 
surface which is normal. Let $p : X \to C$ be an $\A^1$-fibration, where $C$ is an algebraic curve which is either 
affine or projective and $p$ is surjective. We say that the $\A^1$-fibration $p$ is of {\em affine type} (resp. 
{\em complete type}) if $C$ is affine (resp. complete). The $\A^1$-fibration on $X$ is the quotient morphism 
of a $G_a$-action on $X$ if and only if it is of affine type (see \cite{GKM}). We consider the following result 
on deformations. For the complex analytic case, one can refer to \cite{K1} and also to \cite[p. 269]{Iitaka}.

\begin{lem}\label{Lemma 2.1}\label{Lemma 2.1}
Let $\ol{f} : \ol{Y} \to T$ be a smooth projective morphism from a smooth algebraic threefold $\ol{Y}$ to a smooth algebraic 
curve $T$. Let $C$ be a smooth rational complete curve contained in $\ol{Y}_0=\ol{f}^{-1}(t_0)$ for a closed point $t_0$ of 
$T$ \footnote{When we write $t \in T$, we tacitly assume that $t$ is a closed point of $T$}. Then the following assertions hold.
\begin{enumerate}
\item[(1)]
The Hilbert scheme $\Hilb(\ol{Y})$ has dimension less than or equal to $h^0(C, \N_{C/\ol{Y}})$ in the point $[C]$. If 
$h^1(C,\N_{C/\ol{Y}})=0$ then the equality holds and $\Hilb(\ol{Y})$ is smooth at $[C]$. Here $\N_{C/\ol{Y}}$ denotes the 
normal bundle of $C$ in $\ol{Y}$.
\item[(2)]
Let $n=(C^2)$ on $\ol{Y}_0$. Then $\N_{C/\ol{Y}}\cong \SO_C\oplus\SO_C(n)$ provided $n \ge -1$.
\item[(3)]
Suppose $n=0$. Then there exists an \'etale finite covering $\sigma_2 : T' \to T$ such that the morphism $\ol{f}_{T'}$ splits as 
\[
\ol{f}_{T'} : \ol{Y}\times_TT' \st{\varphi} V \st{\sigma_1} T'~,
\]
where $\varphi$ is a $\BP^1$-fibration with $C$ contained as a fiber and $\sigma_1$ makes $V$ a smooth $T'$-scheme of relative 
dimension one with irreducible fibers. Assume further that every smooth rational complete curve $C'$ in $\ol{Y}_0$ satisfies 
$(C'\cdot C)=0$ provided $C'$ is algebraically equivalent to $C$ in $\ol{Y}$. Then the covering $\sigma_2 : T' \to T$ is trivial, 
i.e., $\sigma_2$ is the identity morphism.
\item[(4)]
Suppose $n=-1$. Then $C$ does not deform in the fiber $\ol{Y}_0$ but deforms along the morphism $\ol{f}$ after an \'etale 
finite base change. Namely, there are an \'etale finite morphism $\sigma : T' \to T$ and an irreducible subvariety $Z$ 
of codimension one in $\ol{Y}':=\ol{Y}\times_TT'$ such that $Z$ can be contracted along the fibers of $\ol{f}' : \ol{Y}' 
\to T'$, where $Y'$ is an irreducible smooth affine curve and $\ol{f}'$ is the second projection of $Y\times_TT'$ to $T'$. 
\item[(5)]
Assume that there are no $(-1)$-curves $E$ and $E'$ in $\ol{Y}_0$ such that $E\cap E' \ne \emptyset$ and $E$ is algebraically 
equivalent to $E'$ as $1$-cycles on $\ol{Y}$. Then, after shrinking $T$ to a smaller open set if necessary, we can take 
$Z$ in the assertion (4) above as a subvariety of $\ol{Y}$. The contraction of $Z$ gives a factorization $\ol{f}|_Z : 
Z \st{g} T' \st{\sigma} T$, where $g$ is a $\BP^1$-fibration, $C$ is a fiber of $g$ and $\sigma$ is as above.
\end{enumerate}
\end{lem}

\Proof
(1)~The assertion follows from Grothendieck \cite[Cor. 5.4]{G}.

(2)~We have an exact sequence
\[
0 \to \N_{C/\ol{Y}_0} \to \N_{C/\ol{Y}} \to \N_{\ol{Y}_0/\ol{Y}}|_C \to 0~,
\]
where $\N_{C/\ol{Y}_0}\cong \SO_C(n)$ and $\N_{\ol{Y}_0/\ol{Y}}|_C \cong \SO_C$. The obstruction for this exact sequence 
to split lies in $\Ext^1(\SO_C,\SO_C(n))\cong \H^1(C,\SO_C(n))$, which is zero if $n \ge -1$. 

(3)~Suppose $n=0$. Then $\dim_{[C]}\Hilb(\ol{Y})=2$ and $[C]$ is a smooth point of $\Hilb(\ol{Y})$. Let $H$ be a relatively 
ample divisor on $\ol{Y}/T$ and set $P(n):=P_C(n)=h^0(C,\SO_C(nH))$ the Hilbert polynomial in $n$ of $C$ with respect to $H$. 
Then $\Hilb^P(\ol{Y})$ is a scheme which is projective over $T$. Let $V$ be the irreducible component of $\Hilb^P(\ol{Y})$ 
containing the point $[C]$. Then $V$ is a $T$-scheme with a morphism $\sigma : V \to T$, $\dim V=2$ and $V$ has relative 
dimension one over $T$. Furthermore, there exists a subvariety  $Z$ of $\ol{Y}\times_TV$ such that the fibers of the composite 
morphism 
\[
g : Z \hookrightarrow \ol{Y}\times_TV \st{p_2} V
\]
are curves on $\ol{Y}$ parametrized by $V$. For a general point $v \in V$, the corresponding curve $C':=C_v$ is a smooth 
rational complete curve because $P_{C'}(n)=P(n)$ and $(C')^2=0$ on $\ol{Y}_t=\ol{f}^{-1}(t)$ with $t=\sigma(v)$ because 
$\dim_{[C']}\Hilb(\ol{Y})=2$. In fact, if $(C')^2 \le -1$ then the exact sequence of normal bundles in (2) implies 
$h^0(C',\N_{C'/\ol{Y}}) \le 1$, which contradicts $\dim_{[C']}\Hilb(\ol{Y})=2$.  If $(C')^2 > 0$ then 
$\dim_{[C']}\Hilb(\ol{Y}) >2$, which is again a contradiction. So, $(C')^2=0$. Hence $\ol{Y}_t$ has a $\BP^1$-fibration 
$\varphi_t : \ol{Y}_t \to \ol{B}_t$ such that $C'$ is a fiber and $\ol{B}_t$ is a smooth complete curve. By the 
universality of the Hilbert scheme, there are an open set $U$ of $\ol{B}_t$ and a morphism $\alpha_t : U \to V_t$ such that 
$\varphi_t^{-1}(U)=Z\times_VU$. Since $V$ is smooth over $T$, $\alpha_t$ induces an isomorphism from $\ol{B}_t$ to a connected 
component of $V_t:=\sigma^{-1}(t)$. This is the case if we take $v \in V$ from a different connected component of $V_t$. Let 
$\sigma : V \st{\sigma_1} T' \st{\sigma_2} T$ be the Stein factorization of $\sigma$. Then $\sigma_2$ is an \'etale finite 
morphism and $\sigma_1 : V \to T'$ is a smooth morphism of relative dimension one with irreducible fibers. Furthermore, the 
morphism $g$ above factors as a composite of $T'$-morphsims
\[
g : Z \hookrightarrow (\ol{Y}\times_TT')\times_{T'}V \st{p_2} V \ ,
\]
where $Z$ is identified with $\ol{Y}\times_TT'$ by the above construction. Hence $g$ induces a $T'$-morphism 
$\varphi : \ol{Y}\times_TT' \to V$ such that the composite $\sigma_1\cdot\varphi : \ol{Y}\times_TT' \to V \to T'$ is 
the pull-back $\ol{f}_{T'} : \ol{Y}\times_TT' \to T'$ of the morphism $\ol{f}$. 

In the above argument, we take two curves $C, C'$ corresponding two points $v, v'$ in $V_{t_0}$. Then $C$ is algebraically 
equivalent to $C'$ in $\ol{Y}$, and hence $(C\cdot C')=0$ by the hypothesis. So, $C=C'$ or $C\cap C'=\emptyset$. This 
implies that $C'$ and $C$ are the fibers of the same $\BP^1$-fibration $\varphi_{t_0} : \ol{Y}_{t_0} \to \ol{B}_{t_0}$ 
and hence that $V_{t_0}$ is irreducible. Namely, $\sigma_2^{-1}(t_0)$ consists of a single point. Hence $\deg \sigma_2=1$, 
i.e., $\sigma_2$ is the identity morphism. 

(4)~Suppose $n=-1$. Then $h^0(C,\N_{C/\ol{Y}})=1$ and $h^1(C,\N_{C/\ol{Y}})=0$. Hence $\Hilb^P(\ol{Y})$ has dimension one 
and is smooth at $[C]$, where $P(n)=P_C(n)$ is the Hilbert polynomial of $C$ with respect to $H$. Let $T'$ be the 
irreducible component of $\Hilb^P(\ol{Y})$ containing $[C]$. Note that $\dim T'=1$. Then we find a subvariety $Z$ in 
$\ol{Y}\times_TT'$ such that $C$ is a fiber of $g$ and every fiber of the $T$-morphism $g=p_2|_{T'} : Z \to T'$ is a 
$(-1)$ curve in the fiber $\ol{Y}_t$. In fact, the nearby fibers of $C$ are $(-1)$ curves as a small deformation of $C$ 
by \cite{K1}. Hence, by covering $T'$ by small disks, we know that every fiber of $g$ is a $(-1)$ curve. Further, 
the projection $\sigma : T' \to T$ is a finite morphism as it is projective and $T'$ is smooth because each fiber is 
a $(-1)$ curve in $\ol{Y}$ (see the above arugument for $[C]$). Furthermore, $\sigma$ is \'etale since $\ol{f}$ is locally 
a product of the fiber and the base in the Euclidean topology. Hence $\sigma$ induces a local isomorphism between $T'$ 
and $T$. This implies that $\ol{Y}\times_TT'$ is a smooth affine threefold and the second projection 
$\ol{f}' : \ol{Y}\times_TT' \to T'$ is a smooth projective morphism. Now, after an \'etale finite base change 
$\sigma : T' \to T$, we may assume that $Z$ is identified with a subvariety of $\ol{Y}$. Since $C$ is a $(-1)$ curve 
in $\ol{Y}_0$, it is an extremal ray in the cone $\ol{NE}(\ol{Y}_0)$. Since $C$ is algebraically equivalent to the 
fibers of $g : Z \to T'$, it follows that $C$ is an extremal ray in the relative cone $\ol{NE}(\ol{Y}/T)$. Then 
it follows from \cite[Theorem 3.25]{KollarMori} that $Z$ is contracted along the fibers of $g$ in $\ol{Y}$ and 
the threefold obtained by the contraction is smooth and projective over $T$. 

(5)~Let $\sigma^{-1}(t_0)=\{u_1,\ldots,u_d\}$ and let $Z_{u_i}=Z\cdot(\ol{Y}\times\{u_i\})$ for $1 \le i \le d$. Then 
the $Z_{u_i}$ are the $(-1)$ curves on $\ol{Y}_0$ which are algebraically equivalent to each other as $1$-cycles on $\ol{Y}$. 
By the assumption, $Z_{u_i}\cap Z_{u_j}=\emptyset$ whenever $i \ne j$. This property holds for all $t \in T$ if one 
shrinks to a smaller open set of $t_0$. Then we can identify $Z$ with a closed subvariety of $\ol{Y}$. In fact, the projection 
$p : Z \hookrightarrow Y\times_YT' \to Y$ is a $T$-morphism. For the point $t_0 \in T$, the morphism 
$p\otimes_{\SO_{T,t_0}}\wh{\SO}_{T,t_0}$ with the completion $\wh{\SO}_{T,t_0}$ of $\SO_{T,t_0}$ is a direst sum of 
the closed immersions from $Z\otimes_{\SO_{T,t_0}}\wh{\SO}_{T',u_i}$ into $Y\otimes_{\SO_{T,t_0}}\wh{\SO}_{T',u_i}$ for 
$1 \le i \le r$. So, $p\otimes_{\SO_{T,t_0}}\wh{\SO}_{T,t_0}$ is a closed immersion. Hence $p$ is a closed immersion locally 
over $t_0$ because $\wh{\SO}_{T,t_0}$ is faithfully flat over $\SO_{T,t_0}$. The rest is the same as in the proof of 
the assertion (4).
\QED

Let $Y_0$ be a smooth affine surface and let $\ol{Y}_0$ be a smooth projective surface containing $Y_0$ as an open set in 
such a way that the complement $\ol{Y}_0\setminus Y_0$ supports a reduced effective divisor $D_0$ with simple normal crossings. 
We call $\ol{Y}_0$ a {\em normal completion} of $Y_0$ and $D_0$ the boundary divisor of $Y_0$. An irreducible component of $D_0$ 
is called a {\em $(-1)$ component} if it is a smooth rational curve with self-intersection number $-1$. We say that $\ol{Y}_0$ 
is a {\em minimal} normal completion if the contraction of a $(-1)$ component of $D_0$ (if any) results the image of $D_0$ 
losing the condition of simple normal crossings. 

Let $\ol{f} : \ol{Y} \to T$ be a smooth projective morphism from a smooth algebraic threefold $\ol{Y}$ to a smooth algebraic 
curve $T$ and let $S=\sum_{i=1}^rS_i$ be a reduced effective divisor on $\ol{Y}$ with simple normal crossings. Let $Y= 
\ol{Y}\setminus S$ and let $f =\ol{f}|_Y$. We assume that for every point $t \in T$, the intersection cycle 
$D_t=\ol{f}^{-1}(t)\cdot S$ is a reduced effective divisor of $\ol{Y}_t=\ol{f}^{-1}(t)$ with simple normal crossings 
\footnote{In order to avoid the misreading, it is better to specify our definition of simple normal crossings in the case 
of dimension three. We assume that every irreducible component $S_i$ of $S$ and every fiber $\ol{Y}_t$ are smooth and that 
analytic-locally at every intersection point $P$ of $S_i\cap S_j$ (resp. $S_i\cap S_j\cap S_k$ or $S_i\cap \ol{Y}_t$), $S_i$ 
and $S_j$ (resp. $S_i, S_j$ and $S_k$, or $S_i$ and $\ol{Y}_t$) behave like coordinate hypersurfaces. Hence $S_i\cap S_j$ or 
$S_i\cap \ol{Y}_t$ are smooth curves at the point $P$.} and $Y_t=Y\cap \ol{Y}_t$ is an affine open set of $\ol{Y}_t$. 
For a point $t_0 \in T$, we assume that $\ol{Y}_{t_0}=\ol{Y}_0, D_{t_0}=D_0$ and $Y_{t_0}=Y_0$. A collection 
$(Y,\ol{Y},S,\ol{f},t_0)$ is called a {\em family of logarithmic deformations} of a triple $(Y_0,\ol{Y}_0,D_0)$. We call it 
simply a {\em log deformation} of the triple $(Y_0,\ol{Y}_0,D_0)$. Since $f$ is smooth and $S$ is a divisor with simple normal 
crossings, $(Y,\ol{Y},S,\ol{f},t_0)$ is a family of logarithmic deformations in the sense of Kawamata \cite{Kawamata1, Kawamata}. 

From time to time, we have to make a base change by an \'etale finite morphism $\sigma : T' \to T$ with irreducible $T'$. 
Let $\ol{Y}'=\ol{Y}\times_TT', \ol{f}'=\ol{f}\times_TT', S'=S\times_TT'$ and $Y'=Y\times_TT'$. Since the field extension 
$k(\ol{Y})/k(T)$ is a regular extension, $\ol{Y}'$ is an irreducible smooth projective threefold, and $S'$ is a divisor 
with simple normal crossings. Hence $(Y',\ol{Y}', S', \ol{f}', t_0')$ is a family of logarithmic deformations of the 
triple $(Y'_0,\ol{Y}'_0,D'_0)\cong (Y_0,\ol{Y}_0,D_0)$, where $t_0' \in T'$ with $\sigma(t_0')=t_0$. 
\svskip

We have the following result on logarithmic deformations of affine surfaces with $\A^1$-fibrations.

\begin{lem}\label{Lemma 2.2}
Let $(Y,\ol{Y},S,\ol{f},t_0)$ be a log deformation of the triple $(Y_0,$ $\ol{Y}_0,D_0)$. Then the following assertions hold. 
\begin{enumerate}
\item[(1)]
Assume that $Y_0$ has an $\A^1$-fibration. Then $Y_t$ has an $\A^1$-fibration for every $t \in T$. 
\item[(2)]
If $Y_0$ has an $\A^1$-fibration of affine type (resp. of complete type), then $Y_t$ has also an $\A^1$-fibration of affine 
type (resp. of complete type) for every $t \in T$.
\end{enumerate}
\end{lem}

\Proof
(1)~Note that $K_{\ol{Y}_t}= (K_{\ol{Y}}+\ol{Y}_t)\cdot\ol{Y}_t=K_{\ol{Y}}\cdot\ol{Y}_t$ because $\ol{Y}_t$ is algebraically 
equivalent to $\ol{Y}_{t'}$ for $t' \ne t$. Then $K_{\ol{Y}_t}+D_t=(K_{\ol{Y}}+S)\cdot\ol{Y}_t$. By the hypothesis, 
$h^0(\ol{Y}_0,\SO(n(K_{\ol{Y}_0}+D_0)))=0$ for every $n > 0$. Then the semicontinuity theorem \cite[Theorem 12.8]{H} 
implies that $h^0(\ol{Y}_t,\SO(n(K_{\ol{Y}_t}+D_t)))=0$ for every $n > 0$. Hence $\lkd(Y_t)=-\infty$. Since $Y_t$ is affine, 
this implies that $Y_t$ has an $\A^1$-fibration. 

(2)~Suppose that $Y_0$ has an $\A^1$-fibration $\rho_0 : Y_0 \to B_0$ which is of affine type. Then $\rho_0$ defines a pencil 
$\Lambda_0$ on $\ol{Y}_0$. 

Suppose first that $\Lambda_0$ has no base points and hence defines a $\BP^1$-fibration 
$\ol{\rho}_0 : \ol{Y}_0 \to \ol{B}_0$ such that $\ol{\rho}_0|_{Y_0}=\rho_0$ and $\ol{B}_0$ is a smooth completion of 
$B_0$. If $\ol{\rho}_0$ is not minimal, let $E$ be a $(-1)$ curve contained in a fiber of $\ol{\rho}_0$, which is necessarily 
not contained in $Y_0$. By Lemma \ref{Lemma 2.1}, $E$ extends along the morphism $\ol{f}$ if one replaces the base $T$ by a 
suitable \'etale finite covering $T'$ and can be contracted simultaneously with other $(-1)$ curves contained in the fibers 
$\ol{Y}_t~(t \in T)$. Note that this \'etale finite change of the base curve does not affect the properties of the fiber 
surfaces. Hence we may assume that all simultaneous blowing-ups and contractions as applied below are achieved over the 
base $T$. 

The contraction is performed either within the boundary divisor $S$ or the {\em simultaneous half-point detachments} 
in the respective fibers $Y_t$ for $t \in T$. (For the definition of half-point detachment (resp. attachment), see for example 
\cite{Fujita}). Hence the contraction does not change the hypothesis on the simple normal crossing of $S$ and the intersection 
divisor $S\cdot\ol{Y}_t$. Thus we may assume that $\ol{\rho}_0$ is minimal. Since $B_0 \subsetneqq \ol{B}_0$, a fiber of 
$\ol{\rho}_0$ is contained in a boundary component, say $S_1$. Then the intersection $S_1\cdot \ol{Y}_0$ as a cycle is a disjoint 
sum of the fibers of $\ol{\rho}_0$ with multiplicity one. Hence $(S_1^2\cdot\ol{Y}_0)=((S_1\cdot\ol{Y}_0)^2)_{\ol{Y}_0}=0$. 
Since $\ol{Y}_t$ and $\ol{Y}_0$ are algebraically equivalent, we have $(S_1^2\cdot\ol{Y}_t)=0$ for every $t \in T$. Note that 
$\ol{Y}_t$ is also a ruled surface by Iitaka \cite{Iitaka} and minimal by the same reason as for $\ol{Y}_0$. Considering the 
deformations of a fiber of $\ol{\rho}_0$ appearing in $S_1\cdot\ol{Y}_0$, we know by Lemma \ref{Lemma 2.1} that $S_1\cdot\ol{Y}_t$ 
is a disjoint sum of smooth rational curves with self-intersection number zero. Namely, $S_1\cdot\ol{Y}_t$ is a sum of 
the fibers of a $\BP^1$-fibration. Here we may have to replace the $\BP^1$-fibration $\ol{\rho}_t$ by the second one if 
$\ol{Y}_t \cong \BP^1\times\BP^1$. In fact, if a smooth complete surface has two different $\BP^1$-fibrations and is minimal 
with respect to one fibration, then the surface is isomorphic to $\BP^1\times \BP^1$ and two $\BP^1$-fibrations are the vertical 
and horizontal fibrations. This implies that $Y_t$ has an $\A^1$-fibration of affine type. 

Suppose next that $\Lambda_0$ has a base point, say $P_0$, and that the $\A^1$-fibration $\rho_0$ is of affine type. Then 
all irreducible components of $D_0:=S\cdot\ol{Y}_0$ are contained in the members of $\Lambda_0$. Since the boundary divisor 
$D_0$ of $\ol{Y}_0$ is assumed to be a connected divisor with simple normal crossings, there are at most two components 
of $S\cdot\ol{Y}_0$ passing through $P_0$, and if there are two of them, they lie on different components of $S$ and $P_0$ 
lies on their intersection curve. In particular, if $S_1$ is a component of $S$ containing $P_0$, then $S_1\cdot\ol{Y}_0$ 
is a disjoint sum of smooth rational curves.  Let $C_1$ be the component of $S_1\cdot\ol{Y}_0$ passing through $P_0$ and let 
$F_0$ be the member of $\Lambda_0$ which contains $C_1$. We may assume that $F_0$ is supported by the boundary divisor $D_0$. 
If $F_0$ contains a $(-1)$ curve $E$ such that $P_0 \not\in E$, then $E$ extends along the morphism  $\ol{f}$ and can be 
contracted simultaneously along $\ol{f}$ after the base change by an \'etale finite covering $T' \to T$. So, we may assume 
that every irreducible component of $F_0$ not passing $P_0$ has self-intersection number $\le -2$ on $\ol{Y}_0$. Then we may 
assume that $(C_1^2)_{\ol{Y}_0} \ge 0$. In fact, if there are two irreducible components of $S\cdot\ol{Y}_0$ passing through 
$P_0$ and belonging to the same member $F_0$ of $\Lambda_0$, one of them must have self-intersection number $\ge 0$, 
for otherwise all the components of the member of $\Lambda_0$, after the elimination of base points, would have self-intersection 
number $\le -2$, which is a contradiction. So, we may assume 
that the one on $S_1$, i.e., $C_1$, has self-intersection number $\ge 0$. Then the proper transform of $C_1$ is the unique 
$(-1)$ curve with multiplicity $>1$ in the fiber corresponding to $F_0$ after the elimination of base points of $\Lambda_0$. 

On the other hand, $S_1\cdot\ol{Y}_0$ (as well as $S_i\cdot\ol{Y}_0$ if it is non-empty) is a disjoint sum of smooth rational 
curves, one of which is $C_1$. Let $n:=(C_1^2)_{\ol{Y}_0} \ge 0$. Then $\Hilb^P(\ol{Y})$ has dimension $n+2$ and is smooth 
at the point $[C_1]$. Since $C_1 \cong \BP^1$ and $\N_{C_1/\ol{Y}}\cong \SO(n)\oplus\SO$, $C_1$ extends along the morphism 
$\ol{f}$. Namely, $\ol{f}|_{S_1} : S_1 \to T$ is a composite of a $\BP^1$-fibration $\sigma_1 : S_1 \to T'$ and an \'etale 
finite morphism $\sigma_2 : T' \to T$, where $C_1$ is a fiber of $\sigma_1$. By the base change by $\sigma_2$, we may assume 
that $S_1\cdot\ol{Y}_0=C_1$. In particular, $(C_1^2)_{S_1}=0$. 

Suppose that $C_2$ is a component of $F_0$ meeting $C_1$. Then $C_2$ is contained in a different boundary component, say $S_2$, 
which intersects $S_1$. Since $(H\cdot S_2\cdot \ol{Y}_0) > 0$, we have $(H\cdot S_2\cdot \ol{Y}_t) > 0$ for every $t \in T$, 
where $H$ is a relatively ample divisor on $\ol{Y}$ over $T$. Furthermore, $S_2\cdot \ol{Y}_0$ is algebraically equivalent to 
$S_2\cdot \ol{Y}_t$. Note that $S_2\cdot\ol{Y}_0$ is a disjoint sum of smooth rational curves, one of which is the curve $C_2$ 
connected to $C_1$. By considering the factorization of $\ol{f}|_{S_2} : S_2 \to T$ into a product of a $\BP^1$-fibration and 
an \'etale finite morphism as in the case for $S_1\cdot\ol{Y}_0$ and taking the base change by an \'etale finite morphism, 
we may assume that $S_2\cdot\ol{Y}_0=C_2$. Hence we have 
\[
(S_1\cdot S_2\cdot \ol{Y}_t)=(S_1\cdot S_2\cdot\ol{Y}_0)= (C_1\cdot (S_2\cdot \ol{Y}_0))_{\ol{Y}_0}=(C_1\cdot C_2)_{\ol{Y}_0}=1.
\] 
This implies that $S_2\cdot\ol{Y}_t$ is irreducible for a general point $t \in T$. For otherwise, by the Stein factorization 
of the morphism $\ol{f}|_{S_2} : S_2 \to T$, the fiber $S_2\cdot\ol{Y}_t$ is a disjoint sum $A_1+\cdots+A_s$ of distinct 
irreducible curves which are algebraically equivalent to each other on $S_2$. Since 
\begin{eqnarray*}
\lefteqn{1=(S_1\cdot S_2\cdot \ol{Y}_t) =((S_1\cdot S_2)\cdot (S_2\cdot \ol{Y}_t))_{S_2}} \\
&& =((S_1\cdot S_2)\cdot(A_1+\cdots+A_s))_{S_2}=s((S_1\cdot S_2)\cdot A_1), 
\end{eqnarray*}
we have $s=1$ and $(S_1\cdot S_2\cdot A_1)=1$. So, $\ol{f}|_{S_2} : S_2 \to T$ is now a $\BP^1$-bundle and $(C_2^2)_{S_2}=0$. 
This implies that $\N_{C_2/\ol{Y}}\cong \SO(m)\oplus\SO$ with $m=(C_2^2)_{\ol{Y}_0} \le -2$ and that $C_2$ extends along the 
morphism $\ol{f}$. We can argue in the same way as above with irreducible components of $F_0$ other than $C_1$. 

Assume that no members of $\Lambda_0$ except $F_0$ have irreducible components outside of $Y_0$. If $C_i$ is shown to move on 
the component $S_i$ along the morphism $\ol{f}$, we consider a component $C_{i+1}$ anew which meets $C_i$. Each of them is 
contained in a distinct irreducible boundary component of $S$ and extends along the morphism $\ol{f}$. Let $S_1, S_2, \ldots,S_r$ 
be all the boundary components which meet $\ol{Y}_0$ along the irreducible components of $F_0$. Then $\ol{Y}_t$ intersects 
$S_1+S_2+ \cdots+S_r$ in an effective divisor which has the same form as $F_0$. Furthermore, we have 
\[
((S_i\cdot\ol{Y}_t)^2)_{\ol{Y}_t}=(S_i^2\cdot\ol{Y}_t)=(S_i^2\cdot\ol{Y}_0)=((S_i\cdot\ol{Y}_0)^2)_{\ol{Y}_0}
\] 
for $1 \le i \le r$. Namely, the components $S_i\cdot\ol{Y}_t~(1 \le i \le r)$ with the same multiplicities as $S_i\cdot\ol{Y}_0$ 
in $F_0$ is a member $F_t$ of the pencil $\Lambda_t$ lying outside of $Y_t$. This implies that $\Lambda_t$ has a base point $P_t$ 
and at least one member of $\Lambda_t$ lies outside of $Y_t$. So, the $\A^1$-fibration $\rho_t$ on $Y_t$ is of affine type. 

If the pencil $\Lambda_0$ contains two members $F_0, F_0'$ such that the components $C_1, C_1'$ of $F_0, F'_0$ lie outside of 
$Y_0$ and pass through the point $P_0$, we may assume that $F_0$ is supported by the boundary components, while $F'_0$ may not. 
Then no other members of $\Lambda_0$ have irreducible components outside of $Y_0$ because $\ol{Y}_0\setminus Y_0$ is connected. 
We can argue as above to show, after a suitable \'etale finite base change, that the member $F_0$ moves along the morphism 
$\ol{f}$, and further that every boundary component of $F'_0$ moves on a boundary component, say $S'_j$, as a fiber of 
$f|_{S'_j} : S'_j \to T$. Hence the pencil $\Lambda_t$ has the member $F_t$ corresponding to $F_0$ whose all components lie 
outside of $Y_t$ and the member $F'_t$ corresponding to $F'_0$. In fact, the part of $F'_t$ lying outside of $Y_t$ is determined 
as above, but since $Y_t\cap F'_t$ is a disjoint union of the $\A^1$ which correspond to the $(-1)$ components of $F'_t$ 
(the {\em half-point attachments}), the member $F'_t$ is determined up to its weighted graph. This proof also implies that 
if $\rho_0$ is of complete type then $\rho_t$ is of complete type for every $t \in T$.
\QED

\begin{remark}\label{Remark 2.3}{\em
(1)~In the above proof of Lemma \ref{Lemma 2.2}, the case where the pencil $\Lambda_0$ has a base point $P_0$ on one of the 
connected components $S_1\cap\ol{Y}_0$, say $C_1$, there might exist a monodromy on $\ol{Y}$ which transform $\Lambda_0$ to 
a pencil $\Lambda'_0$ on $\ol{Y}_0$ having a base point $P'_0$ on a different connected component $C'_1$ of $S_1\cap\ol{Y}_0$. 
However, we have $(C_1^2)_{\ol{Y}_0} \ge 0$ as shown in the proof, and $(C_1^2)=({C'}_1^2)$. Since $C'_1$ is contained in 
a member of $\Lambda_0$, whence $({C'}_1^2) < 0$. This is a contradiction. So, $S_1\cap\ol{Y}_0$ is irreducible.

(2)~In the step of the above proof of Lemma \ref{Lemma 2.2} where we assume that no members of $\Lambda_0$ except 
$F_0$ have irreducible components outside of $Y_0$, let $P'_t$ be a point on $C_{1,t}:=S_1\cdot\ol{Y}_t$ other than $P_t$ 
which is the base point of the given pencil $\Lambda_t$. Then there is a pencil $\Lambda'_t$ on $\ol{Y}_t$ which is similar 
to $\Lambda_t$. In fact, note first that $\ol{Y}_t$ is a rational surface. Perform the same blowing-ups with centers at $P'_t$ 
and its infinitely near points as those with centers at $P_t$ and its infinitely near points which eliminate the base points 
of $\Lambda_t$. Then we find an effective divisor $\wt{F}'_t$ supported by the proper transforms of 
$S_i\cdot\ol{Y}_t~(1 \le i \le r)$ and the exceptional curves of the blowing-ups such that $\wt{F}'_t$ has the same form and 
multiplicities as the corresponding member $\wt{F}_t$ in the proper transform $\wt{\Lambda}_t$ of $\Lambda_t$ after the 
elimination of base points. Then $(\wt{F}'_t)^2=0$ and hence $\wt{F}'_t$ is a fiber of an $\BP^1$-fibration on the blown-up 
surface of $\ol{Y}_t$. Then the fibers of the $\BP^1$-fibration form the pencil $\Lambda'_t$ on $\ol{Y}_t$ after the reversed 
contractions. In fact, the surface $Y_t=\ol{Y}_t\setminus D_t$ is the affine plane with two systems of coordinate lines given 
as the fibers of $\Lambda_t$ and $\Lambda'_t$. Hence the $\A^1$-fibrations induced by $\Lambda_t$ and $\Lambda'_t$ are 
transformed by an automorphism of $\ol{Y}_t$. 
\QED}
\end{remark}

The following is one of the simplest examples of our situation.

\begin{example}\label{Example 2.4}{\em
Let $C$ be a smooth conic and let $S$ be the subvariety of codimension one in $\BP^2\times C$ defined by 
\[
S=\{(P,Q) \mid P \in L_Q,~Q \in C\},
\]
where $L_Q$ is the tangent line of $C$ at $Q$. Let $Y=(\BP^2\times C)\setminus S$ and let $f : Y \to C$ be the projection onto 
$C$. We set $T=C$ to fit the previous notations. Set $\ol{Y}=\BP^2\times C$. Then $\ol{f} : \ol{Y} \to T$ is the second 
projection and the boudary divisor $S$ is irreducible. For every point $Q \in C$, $Y_Q:=\BP^2\setminus L_Q$ has a linear pencil 
$\Lambda_Q$ generated by $C$ and $2L_Q$, which induces an $\A^1$-fibration of affine type. The restriction $\ol{f}|_S : S \to T$ 
is a $\BP^1$-bundle. Let $C$ be defined by $X_0X_2=X_1^2$ with respect to a system of homogeneous coordinates $(X_0,X_1,X_2)$ 
of $\BP^2$ and let $\eta=(1,t,t^2)$ be the generic point of $C$ with $t$ an inhomogeneous coordinate on $C \cong \BP^1$. 
Then $L_\eta$ is defined by $t^2X_0-2tX_1+X_2=0$. The generic fiber $Y_\eta$ of $f$ has an $\A^1$-fibration induced by the 
linear pencil $\Lambda_\eta$ whose general members are the conics defined by $(X_0X_2-X_1^2)+u(t^2X_0-2tX_1+X_2)^2=0$, where 
$u \in \A^1$. Indeed, the conics are isomorphic to $\BP^1_{k(t)}$ since they have the $k(t)$-rational point $(1,t,t^2)$, and 
$Y_\eta$ is isomorphic to $\A^2_{k(t)}$. This implies that the affine threefold $Y$ itself has an $\A^1$-fibration. \QED}
\end{example}

In the course of the proof of Lemma \ref{Lemma 2.2}, we frequently used the base change by a finite \'etale morphism 
$\sigma : T' \to T$, where $T'$ is taken in such a way  that for every $t \in T$, the points $\sigma^{-1}(t)$ correspond 
bijectively to the connected components of $S_i\cap \ol{Y}_t$, where $S_i$ is an irreducible component of $S$. Suppose 
that $\deg \sigma > 1$. Let $\ol{Y}'=\ol{Y}\times_TT'$ and $\ol{f}'=\ol{f}\times_TT'$. Note that $\ol{Y}'$ is smooth 
because $\sigma$ is \'etale. The morphism $\sigma : T' \to T$ gives the Stein factorization 
$\ol{f}|_{S_i} : S_i\st{\varphi} T' \st{\sigma} T$. Then the subvariety $S_i$ is considered to be a subvariety of $\ol{Y}'$ 
via a closed immersion $(\id_{S_i},\varphi) : S_i \to S_i\times_TT' \hookrightarrow \ol{Y}\times_TT'$. We denote it by 
$S'_i$. Let $t_1, t_2$ be points of $T'$ such that they correspond to the connected components $A, B$ of $S_i\cap \ol{Y}_t$, 
whence $\sigma(t_1)=\sigma(t_2)=t$. Then $A, B$ are the fibers of $S'_i$ over the points $t_1, t_2$ of $T'$. Hence $A$ and 
$B$ are algebraically equivalent in $\ol{Y}$. Since $T'$ is \'etale over $T$, we say more precisely that they are 
{\em \'etale-algebraically equivalent}. We have $(A^2)_{\ol{Y}_t}=(B^2)_{\ol{Y}_t}$. In fact, noting that $\ol{Y}'_{t_1}$ 
and $\ol{Y}'_{t_2}$ are algebraically equivalent in $\ol{Y}'$ and that $\ol{Y}'_{t_1}$ and $\ol{Y}'_{t_2}$ are isomorphic 
to $\ol{Y}_t$, we have 
\begin{eqnarray*}
(A^2)_{\ol{Y}_t}&=&(A^2)_{\ol{Y}'_{t_1}}=(S'_i\cdot S'_i\cdot\ol{Y}'_{t_1}) \\
&=& (S'_i\cdot S'_i\cdot \ol{Y}'_{t_2})=(B^2)_{\ol{Y}'_{t_2}}=(B^2)_{\ol{Y}_t}
\end{eqnarray*}
Let $C$ be an irreducible curve in $\ol{Y}_0\cap S$, say a connected component of $\ol{Y}_0 \cap S_1$ with an irreducible 
component $S_1$ of $S$. We say that $C$ has {\em no monodromy} in $\ol{Y}$ if $\ol{f}|_{S_1} : S_1 \to T$ has no splitting 
$\ol{f}|_{S_1} : S_1 \st{\sigma_1} T' \st{\sigma_2} T$, where $\sigma_2$ is an \'etale finite morphism with 
$\deg \sigma_2 > 1$. Note that, after a suitable \'etale finite base change $\ol{Y}\times_TT'$, this condition is fulfilled. 
Namely, {\em the monodromy is killed}. Concerning the extra hypothesis in Lemma \ref{Lemma 2.1}, (3) and the possibility of 
achieving the contractions over the base curve $T$ in Lemma \ref{Lemma 2.1}, (5), we have the following result.

\begin{lem}\label{Lemma 2.5}
Let $(Y,\ol{Y},S,\ol{f},t_0)$ be a family of logarithmic deformation of the triple $(Y_0,\ol{Y}_0,D_0)$. Assume that 
$Y_0$ has an $\A^1$-fibration of affine type. Let $\Lambda_0$ be the pencil on $\ol{Y}_0$ whose general members are the 
closures of fibers of the $\A^1$-fibration. Suppose that $\Lambda_0$ defines a $\BP^1$-fibration 
$\varphi_0 : \ol{Y}_0 \to \ol{B}_0$. Suppose further that the section of $\varphi_0$ in $S\cap \ol{Y}_0$ has no monodromy 
in $\ol{Y}$. Then the following assertions hold.
\begin{enumerate}
\item[(1)]
If $C$ is a fiber of $\varphi_0$ with $C\cap Y_0 \ne \emptyset$ and $C'$ is a smooth rational complete curve which is 
algebraically equivalent to $C$ in $\ol{Y}$, then $(C\cdot C')=0$. 
\item[(2)]
There are no two $(-1)$ curves $E_1$ and $E_2$ such that they belong to the same connected component of the Hilbert scheme 
$\Hilb(\ol{Y})$, $E_1$ is an irreducible component of a fiber of $\varphi_0$ and $E_1\cap E_2 \ne \emptyset$.
\end{enumerate}
\end{lem}
\Proof
(1)~Let $S_0$ be an irreducible component of $S$ such that $(S_0\cdot F)=1$ for a general fiber $F$ of $\varphi_0$. Then 
$S_0\cap \ol{Y}_0$ contains a cross-section of $\varphi_0$. The assumption on the absence of the monodromy implies that 
$S_0\cap \ol{Y}_0$ is irreducible and is the section of $\varphi_0$. Note that $\varphi_0$ contains a fiber $F_\infty$ at 
infinity which is supported by the intersection of $\ol{Y}_0$ with the boundary divisor $S$ in $\ol{Y}$. Such a fiber 
exists by the assumption that the $\A^1$-fibration on $Y_0$ is of affine type. Since $S_0\cap\ol{Y}_0$ gives the 
cross-section, $F_\infty$ is supported by $S\setminus S_0$. Since $C\cap(S\setminus S_0)=\emptyset$ and $C'$ is algebraically 
equivalent to $C$ in $\ol{Y}$, $C'$ does not meet the components of $S\setminus S_0$. Hence $C'\cap F_\infty=\emptyset$, and 
$C'$ is a component of a fiber of $\varphi_0$. So, $(C\cdot C')=0$.\footnote{We note here that without the condition on the 
absence of the monodromy of the cross-section, the assertion fails to hold. See Example \ref{Example 2.6}.}

(2)~Suppose that such $E_1$ and $E_2$ exist. Since $E_1$ and $E_2$ are algebraically equivalent $1$-cycles on $\ol{Y}$, $E_1$ 
and $E_2$ have the same intersections with subvarieties of codimension one in $\ol{Y}$. We consider possible cases separately.

(i)~Suppose that both $E_1$ and $E_2$ are contained in the fiber at infinity $F_\infty$. Since $E_1\cap E_2 \ne \emptyset$, 
it follows that $F_\infty=E_1+E_2$ with $(E_1\cdot E_2)=1$. If $E_1$ meets the section $S_0\cap\ol{Y}_0$, then $(E_1\cdot S_0)=1$, 
whence $(E_2\cdot S_0)=1$ because $E_1$ and $E_2$ are algebraically equivalent in $\ol{Y}$. This is a contradiction. Hence 
$E_1\cap E_2=\emptyset$. 

(ii)~Suppose that only $E_1$ is contained in the fiber at infinity $F_\infty$. Take a smooth fiber $F_0$ of $\varphi_0$ with 
$F_0\cap Y_0 \ne \emptyset$ and consider a deformation of $F_0$ in $\ol{Y}$. Then there exist an \'etale finite morphism 
$\sigma_2 : T' \to T$ and a decomposition of $\ol{f}_{T'} : \ol{Y}\times_TT' \st{\varphi} V \st{\sigma_1} T'$ such that $F_0$ 
is a fiber of $\varphi$ (see Lemma \ref{Lemma 2.1}, (3)). Let $B$ be an irreducible curve on $V$ such that $\varphi(F_\infty) 
\not\in B$ and let $W=\varphi^{-1}(B)$. Note that $E_1$ and $E_2$ are also algebraically equivalent in $\ol{Y}\times_TT'$. 
Since $(E_1\cdot W)=0$ by the above construction, it follows that $(E_2\cdot W)=0$. This implies that $E_2$ is contained in a 
fiber of $\varphi_0$. Hence $E_1\cap E_2=\emptyset$. 

(iii)~Suppose that $E_1$ and $E_2$ are not contained in the fiber $F_\infty$. Then $E_1$ and $E_2$ are the fiber 
components of $\varphi_0$ because $(E_i\cdot F_\infty)=0$ for $i=1, 2$. If they belong to the same fiber, we obtain a 
contradiction by the same argument as in the case (i). If they belong to different fibers, then $E_1\cap E_2= \emptyset$. 
\QED

The following example, which is due to one of the referees of this article, shows that Lemma \ref{Lemma 2.5}, (1) does not 
hold without the monodromy condition on the section of $\varphi_0$.

\begin{example}\label{Example 2.6}
{\em Let $Q=\BP^1\times \BP^1$ and $T'=\A^1_*$ which is the affine line minus one point and 
hence is the underlying scheme of the multiplicative group $G_m$. We denote by $\ell$ (resp. $M$) a general fiber of the first 
projection $p_1 : Q \to\BP^1$ (resp. the second projection $p_2 : Q \to \BP^1$). Let $x$ (resp. $y$) be an inhomogeneous 
coordinate on the first (resp. the second) factor of $Q$.  Set $\ell_\infty=p_1^{-1}(\infty)$ and $M_\infty=p_2^{-1}(\infty)$. 
We consider an involution $\iota$ on $Q\times T'$ defined by $(x,y,z) \mapsto (y,x, -z)$, where $z$ is a coordinate of $\A^1_*$.
Let $Q'$ be the blowing-up of $Q$ with center $P_\infty:=\ell_\infty\cap M_\infty$ and let $E$ be the exceptional curve. 
Then the involution $\iota$ extends to the threefold $Q'\times T'$ in such a way that $E\times T'$ is stable under $\iota$. 
Let $\ol{Y}$ be the quotient threefold of $Q'\times T'$ by this $\Z_2$-action induced by the involution $\iota$. Since the 
projection $p_2 : Q'\times T' \to T'$ is $\Z_2$-equivariant, it induces a morphism $\ol{f} : \ol{Y} \to T$, where $T=T'\quot\Z_2
\cong \A^1_*$. Let $S_1=((\ell_\infty \cup M_\infty)\times T')\quot \Z_2$, $S_2=(E\times T')\quot \Z_2$ and $S=S_1+S_2$. Further, 
we let $Y=\ol{Y}\setminus S$ and $f : Y \to T$ the restriction of $\ol{f}$ onto $Y$. Then the following assertions hold.
\begin{enumerate}
\item[(1)]
The surfaces $S_1$ and $S_2$ are smooth irreducible surfaces intersecting normally. 
\item[(2)]
Fix a point $t_0 \in T$ and denote the fibers over $t_0$ with the subscript $0$. Then the collection $(Y,\ol{Y}, S,\ol{f},t_0)$ 
is a family of logarithmic deformations of the triple $(Y_0, \ol{Y}_0, D_0)$, where $D_0=S\cdot\ol{Y}_0$. 
\item[(3)]
For every $t \in T$, $S_1\cap\ol{Y}_t$ is a disjoint union of two smooth curves $C_{1t}, C'_{1t}$ and $S_2\cap \ol{Y}_t$ is a 
smooth rational curve $C_{2t}$, where $(C_{1t}\cdot C_{2t})=(C'_{1t}\cdot C_{2t})=1$ and $(C_{1t}^2)=({C'_{1t}}^2)=(C_{2t}^2)=
-1$. In particular, $C_{1t}$ is \'etale-algebraically equivalent to $C'_{1t}$, and hence has a non-trivial monodromy.
\item[(4)]
Each fiber $\ol{Y}_t$ is isomorphic to $Q'$ with $C_{1t}, C'_{1t}$ and $C_{2t}$ identified with the proper transforms of 
$M_\infty, \ell_\infty$ on $Q'$ and $E$. 
\item[(5)]
Let $\varphi_t : \ol{Y}_t \to \BP_1$ be the $\BP^1$-fibration induced by the first projection $p_1 : Q \to \BP^1$. Then a 
general fiber $\ell=p_1^{-1}(x)$ is algebraically equivalent to $M=p_2^{-1}(x)$ for $x \in T$. 
\item[(6)]
For every $t \in T$, the affine surface $Y_t$ has an $\A^1$-fibration of affine type.
\end{enumerate}}
\end{example}
\Proof
(1)~Since $(Q\setminus(\ell_\infty\cup M_\infty)=\Spec k[x,y]$, the quotient threefold $V=(Q\times T')\quot \Z_2$ contains 
an open set $(\A^2\times T')\quot \Z_2$, which has the coordinate ring over $k$ generated by elements $X=x+y,~U=xy,~~Z=z^2$ and 
$W=(x-y)z$. Hence the open set is a hypersurface $W^2=Z(X^2-4U)$. The quotient threefold $V$ has a similar open neighborhood 
of the image of the curve $\{P_\infty\}\times T'$. This can be observed by taking inhomogeneous coordinates $x', y'$ on $Q$ 
such that $x'=1/x$ and $y'=1/y$, where $\ell_\infty\cup M_\infty$ is given by $x'y'=0$. If we put $W'=x'+y',~U'=x'y'$ and 
$W'=(x'-y')z$, the open neighborhood is defined by a similar equation ${W'}^2=Z({X'}^2-4U')$. Then the image of 
$(\ell_\infty\cup M_\infty)\times T'$ is given by $U'=0$. Hence it has an equation ${W'}^2=Z{X'}^2$. So, this is a smooth 
irreducible surface. The curve $E$ has inhomogeneous coordinate $x'/y'$ (or $y'/x'$). Hence $E$ is stable under the involution 
$\iota$. Note that the involution $\iota$ has no fixed point because there are no fixed points on the factor $T'$. The 
surface $S_1$ is simultaneously contracted along $T$, and by the contraction, $\ol{Y}$ becomes a $\BP^2$-bundle and the surface 
$S_2$ becomes an immersed $\BP^1$-bundle. Then the assertion (1) follows easily.

(2)~The threefold $\ol{Y}$ is smooth and $\ol{f}$ is a smooth morphism. In fact, every closed fiber of $f=\ol{f}|_Y : Y \to 
T$ is isomorphic to the affine plane. 

(3)~If $t=z^2$, $C_{1t}$ (resp. $C'_{1t}$) is identified with $M'_\infty$ (resp. $\ell'_\infty$) in $Q'\times\{z\}$ and 
$\ell'_\infty$ (resp. $M'_\infty$) in $Q'\times\{-z\}$ under the identification $\ol{Y}_t \cong Q'\times\{z\} \cong 
Q'\times\{-z\}$, where $\ell'_\infty$ and $M'_\infty$ are the proper transforms of $\ell_\infty$ and $M_\infty$ on $Q'$.  
Now the rest of the assertions are easily verified.
\QED

A sufficient condition on the absence of the monodromy in Lemma \ref{Lemma 2.5} is given by the following result.

\begin{lem}\label{Lemma 2.7}
Let the notations and the assumptions be the same as in Lemma \ref{Lemma 2.5} and its proof. Let $S_0\cap \ol{Y}_0=
C_{01}\cup \cdots\cup C_{0m}$. Suppose that $C_{01}$ is a section of the $\BP^1$-fibration $\varphi_0$. If $(C_{01}^2) 
\ge 0$, then $C_{01}$ has no monodromy in $\ol{Y}$. Namely, $m=1$ and $S_0\cap \ol{Y}_0$ is irreducible.
\end{lem}
\Proof
Suppose that $m > 1$. Note that $C_{02}, \ldots, C_{0m}$ are mutually disjoint and do not meet a general fiber of $\varphi_0$ 
because they lie outside $Y_0$ and a general fiber meets only $C_{01}$ in the boundary at infinity. This implies that 
$C_{02}, \ldots,C_{0m}$ are rational curves and the fiber components of $\varphi_0$. By the remark given before  
Lemma \ref{Lemma 2.5}, we have $(C_{0i}^2)=(C_{01}^2) \ge 0$ for $2 \le i \le m$. Then $C_{02},\ldots,C_{0m}$ are full 
fibers of $\varphi_0$ and hence they meet the section $C_{01}$. This is a contradiction.
\QED

We prove one of our main theorems.

\begin{thm}\label{Theorem 2.8}
Let $f : Y \to T$ be a morphism from a smooth affine threefold onto a smooth curve $T$ with irreducible general fibers. 
Assume that general fibers of $f$ have $\A^1$-fibrations of affine type. Then, after shrinking $T$ if necessary and taking 
an \'etale finite morphism $T' \to T$, the fiber product $Y'=Y\times_TT'$ has an $\A^1$-fibration which factors the morphism 
$f'=f\times_TT'$. Furthermore, suppose that there is a relative normal completion $\ol{f} : \ol{Y} \to T$ of $f : Y \to T$ 
satisfying the following conditions:
\begin{enumerate}
\item[(1)] 
$(Y,\ol{Y},S,\ol{f},t_0)$ with $t_0 \in T$ and $S=\ol{Y}\setminus Y$ is a family of logarithmic deformations of 
$(Y_0, \ol{Y}_0,D_0)$ as above, where $Y_0=f^{-1}(t_0), \ol{Y}_0=\ol{f}^{-1}(t_0)$ and $D_0=S\cdot\ol{Y}_0$.
\item[(2)]
The given $\A^1$-fibration of affine type on each fiber $Y_t$ extends to a $\BP^1$-fibration $\varphi_t : \ol{Y}_t \to 
\ol{B}_t$.
\item[(3)]
A section of $\varphi_0$ in the fiber $\ol{Y}_0$ lying in $D_0$ has no monodromy in $\ol{Y}$.
\end{enumerate}
Then the given morphism $f : Y \to T$ is factored by an $\A^1$-fibration. 
\end{thm}

\Proof
Embed $Y$ into a smooth threefold $\ol{Y}$ in such a way that $f$ extends to a projective morphism $\ol{f} : \ol{Y} \to T$. 
We may assume that the complement $S:=\ol{Y}\setminus Y$ is a reduced divisor with simple normal crossings. Let $S=S_0+S_1+ 
\cdots+S_r$ be the irreducible decomposition of $S$. For a general point $t \in T$, let $Y_t$ be the fiber $f^{-1}(t)$ and 
let $\rho_t : Y_t \to B_t$ be the given $\A^1$-fibration on $Y_t$. By the assumption, $B_t$ is an affine curve. 
We may assume that $Y_t$ is smooth and hence $B_t$ is smooth. Let $\ol{Y}_t$ be the closure of $Y_t$ in $\ol{Y}$ which we 
may assume to be a smooth projective surface with $t$ a general point of $T$. By replacing $T$ by a smaller Zariski open set, 
we may assume that $\ol{f}$ is a smooth morphism and that $S\cdot\ol{Y}_t$ is a divisor with simple normal crossings for 
every $t \in T$. Hence we may assume that the condition (1) above is realized. 

For each $t \in T$, let $\Lambda_t$ be the pencil generated by the closures (in $\ol{Y}_t$) of the fibers of the 
$\A^1$-fibration $\rho_t$. If $\Lambda_t$ has a base point, we can eliminate the base points by simultaneous blowing ups 
on the boundary at infinity after an \'etale finite base change of $T$. In this step, we may have to replace, for some 
$t \in T$, the pencil $\Lambda_t$ by another pencil $\Lambda'_t$ which also induces an $\A^1$-fibration of affine type on 
$Y_t$ (see the proof of Lemma \ref{Lemma 2.2}). So, we may assume that the condition (2) above is also satisfied. 

If $S_0\cap\ol{Y}_0$ contains a section of $\varphi_0$, we may assume by an \'etale finite base change that $S_0\cap\ol{Y}_0$ 
is irreducible (see the remark before Lemma \ref{Lemma 2.5}). So, we may assume that the condition (3) is satisfied as well. 
 
Hence, we may assume from the beginning that three conditions are satisfied. The fibration $\rho_t$ extends to a $\BP^1$-fibration 
$\varphi_t : \ol{Y}_t \to \ol{B}_t$ for every $t \in T$, where $\ol{B}_t$ is a smooth completion of $B_t$. For $t_0 \in T$, 
we consider the fibration $\varphi_0 : \ol{Y}_0 \to \ol{B}_0$. A general fiber of $\varphi_0$ meets one of the irreducible 
components, say $S_0$, of $S$ in one point. Then so does every fiber of $\varphi_0$ because $S_0\cdot\ol{Y}_0$ is an 
irreducible divisor on $\ol{Y}_0$ and the fibers of $\varphi_0$ are algebraically equivalent to each other on $\ol{Y}_0$. 
Hence $S_0\cdot\ol{Y}_0$ is a section. We claim that 
\begin{enumerate}
\item[(1)]
$\ol{Y}_t$ meets the component $S_0$ for every $t \in T$.
\item[(2)]
After possibly switching the $\A^1$-fibrations if some $Y_t$ has two $\A^1$-fibrations, we may assume that for every 
$t \in T$, the fibers of the $\BP^1$-fibration $\varphi_t$ on $\ol{Y}_t$ meet $S_0$ along a curve $\ol{A}_t$ such that 
$\ol{A}_t$ is a cross-section of $\varphi_t$ and hence $\varphi_t$ induces an isomorphism between $\ol{A}_t$ and $\ol{B}_t$. 
\end{enumerate}

In fact, for a relatively ample divisor $H$ of $\ol{Y}$ over $T$, we have $(H\cdot S_0\cdot\ol{Y}_0) > 0$, whence 
$(H\cdot S_0\cdot \ol{Y}_t)>0$ for every $t \in T$ because $\ol{Y}_t$ is algebraically equivalent to $\ol{Y}_0$. This implies 
the assertion (1). To prove the assertion (2), we consider the deformation of a smooth fiber $C$ of $\varphi_0$ in $\ol{Y}_0$. 
Since general fibers $Y_t$ of $f$ have $\A^1$-fibrations of affine type, by Lemma \ref{Lemma 2.1}, (3) and Lemma \ref{Lemma 2.5}, (1), 
there is a $\BP^1$-fibration $\varphi : \ol{Y} \to V$ such that $C$ is a fiber of $\varphi$. Then the restriction $\varphi|_{\ol{Y}_0}$ 
is the $\BP^1$-fibration $\varphi_0$.  For every $t \in T$, the restriction $\varphi|_{\ol{Y}_t}$ is a $\BP^1$-fibration 
on $\ol{Y}_t$. If it is different from $\varphi_t$, we replace $\varphi_t$ by $\varphi|_{\ol{Y}_t}$. Then 
$(S_0\cdot C')=(S_0\cdot C)=1$ for a general fiber $C'$ of $\varphi_t$ because $C'$ is algebraically equivalent to $C$. 
The assertion follows immediately. 

With the notations in the proof of Lemma \ref{Lemma 2.1}, the isomorphisms $\ol{A}_t \st{\sim} V_t:=\sigma^{-1}(t)\cong \ol{B}_t$ 
shows that the morphism 
\[
S_0 \hookrightarrow \ol{Y} \st{\varphi} V \st{\sigma} T
\]
induces a birational $T$-morphism $S_0 \to V$ and $S_0$ is a cross-section of $\varphi$. It is clear that the boundary divisor $S$ 
contains no other components which are horizontal to $\varphi$. Hence $Y$ has an $\A^1$-fibration. 
\QED 

As a consequence of Theorem \ref{Theorem 2.8}, we have the following result. 

\begin{cor}\label{Corollary 2.9}
Let $f : Y \to T$ be a smooth morphism from a smooth affine threefold $Y$ to a smooth affine curve $T$. Assume that $f$ has 
a relative projective completion $\ol{f} : \ol{Y} \to T$ which satisfies the same conditions on the boundary divisor $S$ and 
the intersection of each fiber $\ol{Y}_t$ with $S$ as set in Lemma \ref{Lemma 2.2}. If a fiber $Y_0$ has a $G_a$-action, 
then there exists an \'etale finite morphism $T' \to T$ such that the threefold $Y'=Y\times_TT'$ has a $G_a$-action as a $T'$-scheme.
Furthermore, if the relative completion $\ol{f} : \ol{Y} \to T$ is taken so that the three conditions in Theorem \ref{Theorem 2.8}
are satisfied, the threefold $Y$ itself has a $G_a$-action as a $T$-scheme.
\end{cor}

\Proof
By Lemma \ref{Lemma 2.2}, every fiber $Y_t$ has an $\A^1$-fibration of affine type $\rho_t : Y_t \to B_t$, where $B_t$ is an 
affine curve. As in the proof of Theorem \ref{Theorem 2.8}, we may assume that the three conditions therein are satisfied. By 
the same theorem, $Y$ has an $\A^1$-fibration $\rho : Y \to U$ such that $f$ is factored as 
\[
f : Y \st{\rho} U \st{\sigma} T~,
\]
where $U_t:=\sigma^{-1}(t)\cong B_t$ for every $t \in T$. Then $U$ is an affine scheme after restricting $T$ to a Zariski open set.
Then $Y$ has a $G_a$-action by \cite{GKM}.
\QED

Given a smooth affine morphism $f : Y \to T$ from a smooth algebraic variety $Y$ to a smooth curve $T$ such that every 
{\em closed} fiber is isomorphic to the affine space $\A^n$ of fixed dimension, one can ask if the generic fiber of $f$ is 
isomorphic to $\A^n$ over the function field $k(T)$. If this is the case with $f$, we say that {\em the generic triviality} holds 
for $f$. In the case $n=2$, this holds by the following theorem. If the generic triviality for $n=2$ holds for $f : Y \to T$ in the 
setup of Theorem \ref{Theorem 2.10}, a theorem of Sathaye \cite{S} shows that $f$ is an $\A^2$-bundle in the sense of Zariski topology. 

\begin{thm}\label{Theorem 2.10}
Let $f : Y \to T$ be a smooth morphism from a smooth affine threefold $Y$ to a smooth affine curve $T$. Assume that the fiber 
$Y_t$ is isomorphic to $\A^2$ for every closed point of $T$. Then the generic fiber $Y_\eta$ of $f$ is isomorphic to the 
affine plane over the function field of $T$. Hence $f : Y \to T$ is an $\A^2$-bundle over $T$ after replacing $T$ by an 
open set if necessary.
\end{thm}

Before giving a proof, we prepare two lemmas where an integral $k$-scheme is a reduced and irreducible algebraic $k$-scheme 
and where a separable $K$-form of $\A^2$ over a field $K$ is an algebraic variety $X$ defined over $K$ such that 
$X\otimes_KK'$ is $K'$-isomorphic to $\A^2$ for a separable algebraic extension $K'$ of $K$.

\begin{lem}\label{Lemma 2.11}
Let $p : X \to T$ be a dominant morphism from an integral $k$-scheme $X$ to an integral $k$-scheme $T$. Assume that the 
fiber $X_t$ is an integral $k$-scheme for every closed point $t$ of $T$. Then the generic fiber $X_\eta=X\times_T\Spec k(T)$ 
is geometrically integral $k(T)$-scheme.
\end{lem}

\Proof
We have only to show that the extension of the function fields $k(X)/k(T)$ is a regular extension. Namely, $k(X)/k(T)$ is a 
separable extension, i.e., a separable algebraic extension of a transcendental extension of $k(T)$ and $k(T)$ is 
algebraically closed in $k(X)$. Since the characteristic of $k$ is zero, it suffices to show that $k(T)$ is algebraically 
closed in $k(X)$. Suppose the contrary. Let $K$ be the algebraic closure of $k(T)$ in $k(X)$, which is a finite algebraic 
extension of $k(T)$. Let $T'$ be the normalization of $T$ in $K$. Let $\nu : T' \to T$ be the normalization morphism which 
is a finite morphism. Then $p : X \to T$ splits as $p : X \st{p'} T' \st{\nu} T$, which is the Stein factorization. Then 
the fiber $X_t$ is not irreducible for a general closed point $t \in T$, which is a contradiction to the hypothesis. 
\QED

The following result is due to Kambayashi \cite{Kamb}.

\begin{lem}\label{Lemma 2.12}
Let $X$ be a separable $K$-form of $\A^2$ for a field $K$. Then $X$ is isomorphic to $\A^2$ over $K$.
\end{lem}

The following proof of Theorem \ref{Theorem 2.10} uses a locally nilpotent derivation and hence is of purely algebraic nature.
\svskip

\noindent
{\bf Proof of Theorem \ref{Theorem 2.10}.}~Every closed fiber $Y_t$ has an $\A^1$-fibration of affine type and hence a 
$G_a$-action. By Corollary \ref{Corollary 2.9}, there exists an \'etale finite morphism $T' \to T$ such that $Y'=Y\times_TT'$ 
has a $G_a$-ation as a $T'$-scheme. Suppose that the generic fiber $Y'_{\eta'}$ of $f_{T'} : Y' \to T'$ is isomorphic to 
$\A^2$ over the function field $k(T')$. Since $Y'_{\eta'}=Y_\eta\otimes_{k(T)}k(T')$, it follows by Lemma \ref{Lemma 2.12} 
that $Y_\eta$ is isomorphic to $\A^2$ over $k(T)$. Hence, we may assume from the beginning that $Y$ has a $G_a$-action which 
induces $\A^1$-fibrations on general closed fibers $Y_t$. The $G_a$-action on a $T$-scheme $Y$ is induced by a locally nilpotent 
derivation $\delta$ on the coordinate ring $B$ of $Y$, i.e., $Y=\Spec B$. Let $T=\Spec R$. Here $\delta$ is an $R$-trivial 
derivation on $B$. Let $A$ be the kernel of $\delta$. Since $B$ is a smooth $k$-algebra of dimension $3$, $A$ is a finitely 
generated, normal $k$-algebra of dimension $2$. The derivation $\delta$ induces a locally nilpotent derivation $\delta_t$ on 
$B_t=B\otimes_RR/\gm_t$, where $\gm_t$ is the maximal ideal of $R$ corresponding to a general point $t$ of $T$. We assume that 
$\delta_t \ne 0$. Since $B_t$ is a polynomial $k$-algebra of dimension $2$ by the hypothesis, $A_t:=\Ker \delta_t$ is a 
polynomial ring of dimension $1$. 

\begin{claim}\label{Claim 1}
$A_t=A\otimes_RR/\gm_t$ if $\delta_t$ is nonzero.
\end{claim}
\Proof
Let $\varphi : B \to B[u]$ be the $k$-algebra homomorphism defined by 
\[
\varphi(b)=\sum_{i \ge 0}\frac{1}{i!}\delta^i(b)u^i~.
\]
Then $\Ker\delta=\Ker(\varphi-\id)$. Hence we have an exact sequence of $R$-modules
\[
0 \to A \to B \st{\varphi-\id} B[u]~. 
\]
Let $\SO_t$ be the local ring of $T$ at $t$, i.e., the localization of $R$ with respect to $\gm_t$, and let $\wh{\SO}_t$ be the 
$\gm_t$-adic completion of $\SO_t$. Since $\wh{\SO}_t$ is a flat $R$-module, we have an 
exact sequence
$$
0 \to A\otimes_R\wh{\SO}_t \to B\otimes_R\wh{\SO}_t \to (B\otimes_R\wh{\SO}_t)[u]~. \eqno{(\ast)}
$$
The completion $\wh{\SO}_t$ as a $k$-module decomposes as $\wh{\SO}_t=k\oplus \wh{\gm}_t$, where $\wh{\gm}_t=\gm_t\wh{\SO}_t$, 
the above exact sequence splits as a direct sum of exact sequences of $k$-modules
\begin{eqnarray*}
&&\ \  0 \to A\otimes_Rk \to B\otimes_Rk \to (B\otimes_Rk)[u] \hspace{5mm} \\
&& 0 \to A\otimes_R\wh{\gm}_t \to B\otimes_R\wh{\gm}_t \to (B\otimes_R\wh{\gm}_t)[u]~.
\end{eqnarray*}

The first one is, in fact, equal to 
\[
0 \to A\otimes_RR/\gm_t \to B_t \st{\varphi_t-\id} B_t[u]~,
\]
where $\varphi_t$ is defined by $\delta_t$ in the same way as $\varphi$ by $\delta$. Hence $\Ker\delta_t=A\otimes_RR/\gm_t=A_t$.
\QED

Let $X=\Spec A$ and let $p : X \to T$ be the morphism induced by the inclusion $R \hookrightarrow A$. Thus $f : Y \to T$ splits as 
\[
f : Y \st{q} X \st{p} T~,
\]
where $q$ is the quotient morphism by the induced $G_a$-action on $Y$. 

\begin{claim}\label{Claim 2}
Suppose that $\delta_t \ne 0$ for every $t \in T$. Then $X$ is a smooth surface with $\A^1$-bundle structure over $T$.
\end{claim}
\Proof
Note that $R$ is a Dedekind domain and $A$ is an integral domain. Hence $p$ is a flat morphism. Since $f$ is surjective, 
$p$ is also surjective. Hence $p$ is a faithfully flat morphism. Further, by Claim \ref{Claim 1}, $X_t=\Spec(A\otimes_RR/\gm_t)$ 
is equal to $\Spec A_t$ for every $t$, which is isomorphic to $\A^1$. In fact, the kernel of a non-trivial locally 
nilpotent derivation on a polynomial ring of dimension $2$ is a polynomial ring of dimension $1$. The generic fiber of $p$ 
is geometrically integral by Lemma \ref{Lemma 2.11}. Hence, by \cite[Theorem 2]{KM}, $X$ is an $\A^1$-bundle over $T$. In 
particular, $X$ is smooth. 
\QED

Let $K=k(T)$ be the function field of $T$. The generic fiber $X_K=X\times_T\Spec K$ is geometrically integral as shown in the 
above proof of Claim \ref{Claim 2}. 

\begin{claim}\label{Claim 3}
The generic fiber $Y_K=Y\times_T\Spec K$ is isomorphic to $\A^2_K$.
\end{claim}

\Proof
We consider $q_K : Y_K \to X_K$, where $X_K \cong \A^1_K$. We prove the following two assertions.
\begin{enumerate}
\item[(1)]
For every closed point $x$ of $X_K$, the fiber $Y_K\times_{X_K}\Spec K(x)$ is isomorphic to $\A^1_{K(x)}$.
\item[(2)]
The generic fiber of $q_K$ is geometrically integral. 
\end{enumerate}

Note that $K(x)$ is a finite algebraic extension of $K$. Let $T'$ be the normalization of $T$ in $K':=K(x)$. We consider 
$Y':=Y\times_TT'$ instead of $Y$. Then the $G_a$-action on $Y$ lifts to $Y'$ and the quotient variety is $X'=X\times_TT'$. 
Indeed, the normalization $R'$ of $R$ in $K'$ is the coordinate ring of $T'$ and is a flat $R$-module. Then the sequence 
of $R'$-modules
\[
0 \to A\otimes_RR' \to B\otimes_RR' \st{\varphi'-\id} (B\otimes_RR')[u]
\]
is exact, where $\varphi'=\varphi\otimes_RR'$. Hence $q_{K'} : Y'_{K'} \to X'_{K'}$, which is the base change of $q_K$ 
with respect to the field extension $K'/K$, is the quotient morphism by the $G_a$-action on $Y'_{K'}$ induced by $\delta$. 
Since $X'_{K'}=X\times_T\Spec K'$, there exists a $K'$-rational point $x'$ on $X'_{K'}$ such that $x$ is the image of 
$x'$ by the projection $X'_{K'} \to X_K$. If the fiber of $q_{K'}$ over $x'$, i.e., $Y'_{K'}\times_{X'_{K'}}(\Spec K', x')$,  
is isomorphic to $\A^1_{K'}$, then $Y_K\times_{X_K}\Spec K'$ is isomorphic to $\A^1_{K'}$ because 
$Y'_{K'}\times_{X'_{K'}}\Spec K'=Y_K\times_{X_K}\Spec K'$. Thus we may assume that $x$ is a $K$-rational
point. Let $C$ be the closure of $x$ in $X$. Then $C$ is a cross-section of $p : X \to T$. Let $Z:=Y\times_XC$. Then 
$q_C : Z \to C$ is a faithfully flat morphism such that the fiber $q_C^{-1}(w)$ is isomorphic to $\A^1$ for every closed 
point $w \in C$. In fact, $q_C^{-1}(w)$ is the fiber of $Y_t \to X_t$ over the point $w \in C$, where $t=p(w), Y_t \cong \A^2, 
X_t \cong \A^1$ and $X_t=Y_t\quot G_a$. By Lemma \ref{Lemma 2.11} (which is extended to a non-closed field $K$), 
the generic fiber of $q_C$ is geometrically integral, and the generic fiber of $q_C$, which is $Y_K\times_{X_K}\Spec K(x)$, 
is isomorphic to $\A^1_{K'}$ by \cite[Theorem 2]{KM}. This proves the first assertion.

The generic point of $X_K$ corresponds to the quotient field $L:=Q(A)$. Then it suffices to show that $B\otimes_AQ(A)$ is 
geometrically integral over $Q(A)$. Meanwhile, $B\otimes_AQ(A)$ has a locally nilpotent derivation $\delta\otimes_AQ(A)$ 
such that $\Ker(\delta\otimes_AQ(A))=Q(A)$. Hence $B\otimes_AQ(A)$ is a polynomial ring $Q(A)[u]$ in one variable over 
$Q(A)$ because $\delta\otimes_AQ(A)$ has a slice. So, $B\otimes_AQ(A)$ is geometrically integral over $Q(A)$. Now, by 
\cite[Theorem]{KW}, $Y_K$ is an $\A^1$-bundle over $X_K\cong \A^1_K$. Hence $Y_K$ is isomorphic to $\A^2_K$. We have 
to replace $T$ by an open set $T\setminus F$, where $F=\{t \in T \mid \delta_t=0\}$. This completes the proof of 
Theorem \ref{Theorem 2.10}.
\QED

We can prove Theorem \ref{Theorem 2.10} in a more geometric way by making use of a theorem of Ramanujam-Morrow on the 
boundary divisor of a minimal normal completion of the affine plane \cite{R, Morrow}. The proof given below is explained 
in more precise and explicit terms in \cite[Lemma 3.2]{KZ}. In particular, the step to show that $\ol{Y}_K \cong \BP^2_K$ 
and $Y_K \cong \A^2_K$ is due to [{\em loc.cit.}].
\svskip

\noindent
{\bf The second proof of Theorem \ref{Theorem 2.10}.}~
Let $f:Y \to T$ be as in Theorem \ref{Theorem 2.10}. Let $\ol{Y}$ be a relative completion such that $\ol{Y}$ is smooth 
and $f$ extends to a smooth projective morphism $\ol{f} : \ol{Y} \to T$ with the conditions in Lemma \ref{Lemma 2.2} 
being satisfied together with $S:= \ol{Y}\setminus Y$. To obtain this setting, we may have to shrink $T$ to a smaller 
open set of $T$. As in the first proof and the proof of Lemma \ref{Lemma 2.2}, we can apply an \'etale finite base change 
$T' \to T$ by which the intersection $S_i\cap \ol{Y}_t$ is irreducible for every irreducible component $S_i$ of $S$ and 
every $t \in T$. In particular, we assume 
that $\ol{Y}_t$ is a smooth normal completion of $Y_t$ for every $t \in T$, where $Y_t$ is isomorphic to $\A^2$. Fix one 
such completed fiber, say $\ol{Y}_0=\ol{f}^{-1}(t_0)$, and consider the reduced effective divisor $\ol{Y}_0\setminus Y_0$ 
with $Y_0=f^{-1}(t_0)\cong \A^2$. Namely, $(Y,\ol{Y},S,\ol{f},t_0)$ is a log deformation of $(Y_0,\ol{Y}_0,D_0)$. 
If the dual graph of this divisor is not linear then it contains a $(-1)$-curve meeting at most two other components of 
$D_0$ by a result of Ramanujam \cite{R}. By (4) of Lemma \ref{Lemma 2.1}, such a $(-1)$-curve deforms along the fibers of 
$\ol{f}$ and we get an irreducible component, say $S_1$, of $S=\sum_{i=0}^rS_i$ which can be contracted. Repeating this 
argument, we can assume that all the dual graphs for $\ol{Y}_t\setminus Y_t$, as t varies on the set of closed points of 
$T$, are linear chains of smooth rational curves. By \cite{Morrow}, at least one of these curves is a $(0)$-curve. Fix 
such a $(0)$-curve $C_1$ in $\ol{Y}_0\setminus Y_0$. Then $C_1$ deforms along the fibers of $\ol{f}$ and forms an 
irreducible component, say $S_1$, of $S$ by abuse of the notations. By the argument in the proof of Lemma \ref{Lemma 2.2}, 
if $C_2$ is a component of $\ol{Y}_0\setminus Y_0$ meeting $C_1$, it deforms along the fibers of $\ol{f}$ on an 
irreducible component, say $S_2$, of $S$. Repeating this argument, we know that all irreducible components of 
$\ol{Y}_0\setminus Y_0$ extend along the fibers of $\ol{f}$ to form the irreducible components of $S$ and that the dual 
graphs of $\ol{Y}_t\setminus Y_t$ are the same for every $t \in T$. Now let $K$ be the function field of $T$ over $k$. 
We consider the generic fibers $\ol{Y}_K$ and $Y_K$ of $\ol{f}$ and $f$. Then the dual graph of $\ol{Y}_K\setminus Y_K$ 
is the same linear chain of smooth rational curves as the closed fibers $\ol{Y}_t\setminus Y_t$. Write 
$\ol{Y}_0\setminus Y_0=\sum_{i=1}^rC_i$. If $C_i$ and $C_j$ meet for $i \ne j$, then the intersection point $C_i\cap C_j$ 
moves on the intersection curve $S_i\cdot S_j$. Since any minimal normal completion of $\A^2$ can be brought to $\BP^2$ 
by blowing ups and downs with centers on the boundary divisor, we can blow up simultaneously the intersection curves and 
blow down the proper transforms of the $S_i$ according to the blowing ups and downs on $\ol{Y}_0$. Here we note that 
the begining center of blowing up is a point on a $(0)$-curve $C_1$. In this case, we choose a suitable cross-section 
on the irreducible component $S_1$ which is a $\BP^1$-bundle in the Zariski topology because $\dim T=1$. Note that if 
$T$ is irrational, then the chosen cross-section may meet the intersection curves on $S_1$ with other components of $S$. 
Then we shrink $T$ so that the cross-section does not meet the intersection curves. If $T$ is rational, $S_1$ is a trivial 
$\BP^1$-bundle, hence we do not need the procedure of shrinking $T$. Thus we may assume that, for every $t \in T$, $\ol{Y}_t$ 
is isomorphic to $\BP^2$ and $\ol{Y}_t\setminus Y_t$ is a single curve $C_t$ with $(C_t)^2=1$. This implies that 
$\ol{Y}_K \cong \BP^2_K$ and $Y_K \cong \A^2_K$. 
\QED

In connection with Theorem \ref{Theorem 2.10}, we can pose the following 

\begin{prob}\label{Problem 2.13}
Let $K$ be a field of characteristic zero and let $X$ be a smooth affine surface defined over $K$. Suppose that 
$X\otimes_K\ol{K}$ has an $\A^1$-fibration of affine type, where $\ol{K}$ is an algebraic closure of $K$. Does $X$ then have
an $\A^1$-fibration of affine type?
\end{prob}

If we consider an $\A^1$-fibration of complete type, an example of Dubouloz-Kishimoto gives a counter-example to a  
similar problem for the complete type (see Theorem \ref{Theorem 5.1}). In view of Example \ref{Example 2.6} and Theorem 
\ref{Theorem 2.8}, we need perhaps some condition for a positive answer in the case of affine type which guarantees the 
absence of monodromy of a cross-section of a given $\A^1$-fibration.

\section{Topological arguments instead of Hilbert schemes}
  
In this section we will briefly indicate topological proofs of some of the results in the section two. 
The use of topological arguments would make the cumbersome geometric arguments more transparent for 
the people who do not appreciate the heavy machinery like Hilbert scheme. 

We will use the following basic theorem due to Ehresmann \cite[Chapter V, Prop. 6.4]{W}.

\begin{thm}\label{Theorem 3.1}
Let $M$ be a connected differentiable manifold, $S$ a closed submanifold, $f : M \to N$ a proper differentiable map 
such that the tangent maps corresponding to $f$ and $f|_S : S \to N$ are surjective at any point in $M$ and $S$. 
Then $f|_{M\setminus S} : M\setminus S \to N$ is a locally trivial fiber bundle with respect to the base $N$. 
\end{thm}

Note that the normal bundle of any fiber of $f$ is trivial. We can give a proof of Ehresmann's theorem using this 
observation, and the well-known result from differential topology that given a compact submanifold $S$ of a 
$C^{\infty}$ manifold $X$ there are arbitrarily small tubular neighborhoods of $S$ in $X$ which are diffeomorphic 
to neighborhoods of $S$ in the total space of normal bundle of $S$ in $X$ \cite[Chapter II, Theorem 11.14]{B}.
\svskip

Now let $\ol{f} : \ol{Y} \to T$ be a smooth projective morphism from a smooth algebraic threefold onto a smooth 
algebraic curve $T$. Let $\ol{Y}_t=\ol{f}^{-1}(t)$ be the fiber over $t \in T$. Let $S$ be a simple normal crossing 
divisor on $\ol{Y}$ such that $D_t:= S\cap \ol{Y}_t$ is a simple normal crossing divisor for each $t \in T$ and 
$Y_t:=\ol{Y}_t\setminus D_t$ is affine for each $t$.
\svskip

We can assume that $\ol{f} : \ol{Y} \to T$ has the property that the tangent map is surjective at each point. It 
follows from Ehresmann's theorem that all the surfaces $\ol{Y}_t$ are mutually diffeomorphic. In particular, they have 
the same topological invariants like the fundamental group $\pi_1$ and the Betti number $b_i$. By shrinking $T$ if necessary, 
we will assume that the restricted map $\ol{f} : S_i \to T$ is smooth for each $i$. For fixed $i$ and $t_0$ the intersection 
$S_i\cap \ol{Y}_{t_0}$ is a disjoint union of smooth, compact, irreducible curves. Let $C_{t_0,i}$ be one of these 
irreducible curves. Then for each $t$ which is close to $t_0$, there is an irreducible curve $C_{t,i}$ in 
$S_i\cap \ol{Y}_t$ and suitable tubular neighborhoods of $C_{t_0,i},~C_{t,i}$ in $\ol{Y}_{t_0},~\ol{Y}_{t}$ respectively 
are diffeomorphic by Ehresmann's theorem. This implies that $C_{t_0,i}^2$ in $\ol{Y}_{t_0}$ and $C_{t,i}^2$ in 
$\ol{Y}_t$ are equal. This proves that the weighted dual graphs of the curves $D_t$ in $\ol{Y}_t$ are the same 
for each $t$.
\svskip

Recall that if $X$ is a smooth projective surface with a smooth rational curve $C\subset X$ such that $C^2=0$ 
then $C$ is a fiber of a $\BP^1$-fibration on $X$. If the irregularity $q(X)>0$ then the Albanese morphism 
$X \to \Alb(X)$ gives a $\BP^1$-fibration on $X$ with $C$ as a fiber. By the above discussion the fiber surfaces 
$\ol Y_t$ have the same irregularity.
\svskip

Suppose that $\ol{Y}_0$ has an $\A^1$-fibration of affine type $f : Y_0 \to B$. If $\ol{f} : \ol{Y}_0 \to \ol{B}$ is an 
extension of $f$ to a smooth completion of $Y_0$ then, after simultaneous blowing ups and downs along the fibers of 
$\ol{f}$,  we may asssume that $D_0:=\ol Y_0\setminus Y_0$ contains at least one $(0)$-curve which is a tip, i.e., 
the end component of a maximal twig of $D_0$. Since $D_t$ and $D_0$ have the same weighted dual graphs $D_t$ also 
contains a $(0)$-curve which is a tip of $D_t$. Hence, $Y_t$ also has an $\A^1$-fibration of affine type. This proves 
the assertion (2) in Lemma 2.2.
\svskip 

We can also shorten the part of showing the invariance of the boundary weighted graphs in the second proof of Theorem 
\ref{Theorem 2.10}. Suppose now that $f : Y \to T$ is a fibration on a smooth affine threefold $Y$ onto a smooth curve $T$ 
such that every scheme-theoretic fiber of $f$ is isomorphic to $\A^2$. We can embed $Y$ in a smooth projective threefold 
$\ol{Y}$ such that $f$ extends to a morphism $\ol{f} : \ol{Y} \to T$. By shrinking $T$ we can assume that $\ol{f}$ is smooth, 
each irreducible component $S_i$ of $\ol{Y}\setminus Y$ intersects each $\ol{Y}_t$ transversally, etc. By the above 
discussions, each $D_t:=\ol{Y}_t\setminus Y_t$ has the same weighted dual graph. Since $Y_t$ is isomorphic to $\A^2$, 
we can argue as in the second proof of Theorem \ref{Theorem 2.10} using the result of Ramanujam-Morrow to conclude that 
$f$ is a trivial $\A^2$-bundle on a non-empty Zariski-open subset of $T$. This observation applies also to the proof of 
Theorem \ref{Theorem 4.6}.

\section{Deformations of $\ML_0$ surfaces}

For $i=0, 1, 2$, an $\ML_i$ surface is by definition a smooth affine surface $X$ such that the Makar-Limanov invariant 
$\ML(X)$ has transcendence degree $i$ over $k$ \cite{GMMR}. In this section, we assume that the ground field $k$ is the 
complex field $\C$. Let $\SF=(Y,\ol{Y},S,\ol{f},t_0)$ be a family satisfying the conditions of Lemma \ref{Lemma 2.2}. 
Let $D_0=S\cap \ol{Y}_0$. 

\begin{lem}\label{Lemma 4.1}
Let $\SF=(Y,\ol{Y},S,\ol{f}, t_0)$ be a log deformation of $(Y_0, \ol{Y}_0,D_0)$. Assume that $D_0$ is a tree of smooth 
rational curves satisfying one of the following conditions.
\begin{enumerate}
\item[(i)]
$D_0$ contains an irreducible component $C_1$ such that $(C_1^2) \ge 0$.
\item[(ii)]
$D_0$ contains a $(-1)$ curve which meets more than two other components of $D_0$. 
\end{enumerate}
Then the following assertions hold after changing $T$ by an \'etale finite covering of an open set of $T$ if necessary.
\begin{enumerate}
\item[(1)]
Every irreducible component of $D_0$ deforms along the fibers of $\ol{f}$. Namely, if $D_0=\sum_{i=1}^rC_i$ is the 
irreducible decomposition, then, for every $1 \le i \le r$, there exists an irreducible component $S_i$ of $S$ such 
that $\ol{f}|_{S_i} : S_i \to T$ has the fiber $(\ol{f}|_{S_i})^{-1}(t_0)=C_i$. Furthermore, $S=\sum_{i=1}^r S_i$.
\item[(2)]
For $t \in T$, let $C_{i,t}=(\ol{f}|_{S_i})^{-1}(t)$. Then $D_t=\sum_{i=1}^rC_{i,t}$ and $D_t$ has the same weighted 
graph on $\ol{Y}_t$ as $D_0$ does on $\ol{Y}_0$.
\item[(3)]
For every $i$, $\ol{f}|_{S_i} : S_i \to T$ is a trivial $\BP^1$-bundle over $T$.
\end{enumerate}
\end{lem}
 
\Proof
By a suitable \'etale finite base change of $T$, we may assume that $S_i\cap \ol{Y}_0$ is irreducible for every 
irreducible component $S_i$ of $S$. Then the argument is analytic locally almost the same as in the proof for the 
assertion (2) of Lemma \ref{Lemma 2.2}. Consider the deformation of $C_1$ along the fibers of $\ol{f}$, which moves 
along the fibers because $(C_1^2) \ge -1$. Then the components of $D_0$ which are adjacent to $C_1$ also move along 
the fibers of $\ol{f}$. Once these components of $D_0$ move, then the components adjacent to these components move 
along the fibers of $\ol{f}$. Since $D_0$ is connected because $Y_0$ is affine, all the components of $D_0$ move 
along the fibers of $\ol{f}$. If $S$ contains an irreducible component which does not intersect $\ol{Y}_0$, it is 
a fiber component of $\ol{f}$. Then we remove the fiber by shrinking $T$. This proves the assertion (1). 

Let $S=\sum_{i=1}^rS_i$ be the irreducible decomposition of $S$. As shown in (1), $S_i\cap \ol{Y}_0 \ne \emptyset$ 
for every $i$. Then $S_i\cap \ol{Y}_t \ne \emptyset$ as well by the argument in the proof of Lemma \ref{Lemma 2.2}.

Note that $((S_i\cdot \ol{Y}_t)^2)_{\ol{Y}_t}=(S_i^2\cdot \ol{Y}_t)=(S_i^2\cdot\ol{Y}_0)=((S_i\cdot\ol{Y}_0)^2)_{\ol{Y}_0}$ 
because $\ol{Y}_t$ is algebraically equivalent to $\ol{Y}_0$. Hence $D_0$ and $D_t$ have the same dual graphs.
\QED

In order to prove the following result, we use Ehresmann's theorem, which is Theorem \ref{Theorem 3.1}. 

\begin{lem}\label{Lemma 4.2}
Let $\SF=(Y, \ol{Y},S,\ol{f},t_0)$ be a log deformation of $(Y_0,\ol{Y}_0,D_0)$ which satisfies the same conditions as 
in Lemma \ref{Lemma 4.1}. Assume further that $p_g(\ol{Y}_0)=q(\ol{Y}_0)=0$. Then the following assertions hold.
\begin{enumerate}
\item[(1)]
$\Pic(Y_t)\cong \Pic(Y_0)$ for every $t \in T$.
\item[(2)]
$\Gamma(Y_t,\SO_{Y_t}^*)\cong \Gamma(Y_0,\SO_{Y_0}^*)$ for every $t \in T$.
\end{enumerate}
\end{lem}

\Proof
Since $p_g$ and $q$ are deformation invariants, we have $p_g(\ol{Y}_t)=q(\ol{Y}_t)=0$ for every $t \in T$. The 
exact sequence
\[
0 \lto \Z \lto \SO_{\ol{Y}_t} \st{\exp} \SO_{\ol{Y}_t}^* \lto 0
\]
induces an exact sequence 
\[
H^1(\ol{Y}_t,\SO_{\ol{Y}_t}) \to H^1(\ol{Y}_t,\SO_{\ol{Y_t}}^*) \to H^2(\ol{Y}_t;\Z) \to H^2(\ol{Y}_t,\SO_{\ol{Y}_t})
\]
Since $p_g(\ol{Y}_t)=q(\ol{Y}_t)=0$, we have an isomorphism 
\[
H^1(\ol{Y}_t,\SO_{\ol{Y}_t}^*)\cong H^2(\ol{Y}_t;\Z)\ .
\]
Now consider the canonical homomorphism $\theta_t : H_2(D_t;\Z) \to H_2(\ol{Y}_t;\Z)$, where $H_2(\ol{Y}_t;\Z)\cong 
H^2(\ol{Y}_t;\Z)=\Pic(\ol{Y}_t)$ by the Poincar\'e duality. Then $\coim\theta_t=\Pic(Y_t)$ and $\Ker\theta_t=
\Gamma(Y_t,\SO_{Y_t}^*)/k^*$. 

Let $N$ be a nice tubular neighborhood of $S$ with boundary in $\ol{Y}$. The smooth morphism $\ol{f} : \ol{Y} \to T$ 
together with its restriction on the $(N,\partial N)$ gives a proper differential mapping which is surjective and 
submersive. By Theorem \ref{Theorem 3.1}, it is differentiably a locally trivial fibration. Namely, there exists 
a small disc $U$ of $t_0$ in $T$ and a diffeomorphism $\varphi_0 : \ol{Y}_0\times U \st{\approx} (\ol{f})^{-1}(U)$ 
such that its restriction induces a diffeomorphism 
\[
\varphi_0 : (N\cap\ol{Y}_0)\times U \st{\approx} (\ol{f}|_N)^{-1}(U)\ .
\]
For $t \in U$, noting that $U$ is contractible and hence $H_2(\ol{Y}_0\times U;\Z)=H_2(\ol{Y}_0;\Z)$ and 
$H_2((N\cap\ol{Y}_0)\times U;\Z)=H_2(N\cap\ol{Y}_0;\Z)$, the inclusions $\ol{Y}_t \hookrightarrow (\ol{f})^{-1}(U)$ and 
$N\cap\ol{Y}_0 \hookrightarrow (\ol{f}|_N)^{-1}(U)$ induces compatible isomorphisms 
\[
p_t : H_2(\ol{Y}_t;\Z) \to H_2((\ol{f})^{-1}(U);\Z) \st{(\varphi^{-1})_*} H_2(\ol{Y}_0\times U;\Z)=H_2(\ol{Y}_0;\Z)
\]
and its restriction $q_t : H_2(N\cap \ol{Y}_t;\Z) \st{\sim} H_2(N\cap\ol{Y}_0;\Z)$. Since $S$ and hence $D_t$ are 
strong deformation retracts of $N$ and $N\cap\ol{Y}_t$ respectively, the isomorphism $q_t$ induces an isomorphism 
$r_t : H_2(D_t;\Z) \st{\sim} H_2(D_0;\Z)$ such that the following diagram 
\[
\CD
H_2(D_t;\Z) @>\theta_t>> H_2(\ol{Y}_t;\Z) \\
@Vr_tVV    @VVp_tV \\
H_2(D_0;\Z) @>>\theta_0> H_2(\ol{Y}_0;\Z) 
\endCD
\]
This implies that $\Pic(Y_t) \cong \Pic(Y_0)$ and $\Gamma(Y_t,\SO_{Y_t}^*)\cong \Gamma(Y_0,\SO_{Y_0}^*)$. If $t$ is 
an arbitrary point of $T$, we choose a finite sequence of points $\{t_0, t_1, \ldots, t_n=t\}$ such that $t_i$ is in 
a small disc $U_{i-1}$ around $t_{i-1}\ (1 \le i \le n)$ for which we can apply the above argument. 
\QED

\begin{remark}\label{Remark 4.3}{\em 
By a result of W. Neumann \cite[Theorem 5.1]{N}, if $X$ is a normal affine surface, $D$ an SNC divisor at infinity 
of $X$ which does not contain any $(-1)$-curve meeting at least three other components of $D$ and all whose maximal 
twigs are smooth rational curves with self-intersections $\le -2$, then the boundary $3$-manifold of a nice tubular 
neighborhood $N$ of $D$ determines the dual graph of $D$. If we use the local differentiable triviality of a tubular 
neighborhood $N$, this result of Neumann shows that the weighted dual graph of $D_t$ is deformation invariant.}
\end{remark}

According to \cite[Lemmas 1.2, 1.4]{GMMR}, we have the following property and characterization of $\ML_0$-surface.

\begin{lem}\label{Lemma 4.4}
Let $X$ be a smooth affine surface and let $V$ be a minimal normal completion of $X$. Then the following assertions 
hold.
\begin{enumerate}
\item[(1)]
$X$ is an $\ML_0$-surface if and only if $\Gamma(X,\SO_X^*)=k^*$ and the dual graph of the boundary divisior $D:=
V-X$ is a linear chain of smooth rational curves.
\item[(2)]
If $X$ is an $\ML_0$-surface, $X$ has an $\A^1$-fibration, and any $\A^1$-fibration $\rho : X \to B$ has base curve 
either $B\cong \BP^1$ or $B \cong \A^1$. If $B \cong \BP^1$, $\rho$ has at most two multiple fibers, and if $B 
\cong \A^1$, it has at most one multiple fiber.
\end{enumerate}
\end{lem}

The following result is a direct consequence of the above lemmas.

\begin{thm}\label{Theorem 4.5}
Let $\SF=(Y,\ol{Y},S,\ol{f},t_0)$ be a log deformation of $(Y_0,\ol{Y}_0,D_0)$, where $Y_0$ is an $\ML_0$-surface. 
Then $Y_t$ is an $\ML_0$-surface for every $t \in T$.
\end{thm}

\Proof
If $S\cap \ol{Y}_t$ contains a $(-1)$ curve, then it deforms along the fibers of $\ol{f}$ after an \'etale finite 
base change of $T$, and these $(-1)$ curves are contracted simultaneously by Lemma \ref{Lemma 2.1}. Hence we may 
assume that $\ol{Y}_t$ is a minimal normal completion of $Y_t$ for every $t \in T$. By Lemma \ref{Lemma 4.4}, 
$D_0:=S\cap \ol{Y}_0$ is a linear chain of smooth rational curves. Hence $D_t:=S\cap\ol{Y}_t$ is also a linear 
chain of smooth rational curves. By Lemma \ref{Lemma 4.2}, $\Gamma(Y_t,\SO_{Y_t}^*)=k^*$ for every $t \in T$ 
because $\Gamma(Y_0,\SO_{Y_0}^*)=k^*$. So, $Y_t$ is an $\ML_0$-surface by Lemma \ref{Lemma 4.4}.
\QED 

A smooth affine surface $X$ is, by definition, an {\em affine pseudo-plane} if it has an $\A^1$-fibration of affine type 
$p : X \to \A^1$ admitting at most one multiple fiber of the form $m\A^1$ as a singular fiber (see \cite{M2} for 
the definition and relevant results). An affine pseudo-plane is a $\Q$-homology plane, its Picard group is a cyclic group 
$\Z/m\Z$ and there are no non-constant invertible elements. An $\ML_0$-surface is an affine pseudo-plane if the Picard 
number is zero.

If $\ol{X}$ is a minimal normal completion of an affine pseudo-plane $X$, the boundary divisor $D=\ol{X}-X$ is a tree 
of smooth rational curves, which is not necessarily a linear chain. By blowing-ups and blowing-downs with centers on 
the boundary divisor $D$, we can make the completion $\ol{X}$ satisfy the following conditions \cite[Lemma 1.7]{M2}.
\begin{enumerate}
\item[(i)]
There is a $\BP^1$-fibration $\ol{p} : \ol{X} \to \BP^1$ which extends the $\A^1$-fibration $p : X \to \A^1$.
\item[(ii)]
The weighted dual graph of $D$ is 
\[
\begin{array}{ccccc}
(0) & \raisebox{1mm}{\rule{6mm}{0.2mm}} & (0) & \raisebox{1mm}{\rule{6mm}{0.2mm}} & A \\
\mbox{} \\
\ell & & M 
\end{array}
\]
\item[(iii)]
There is a $(-1)$ curve $F_0$ (called {\em feather}) such that $F_0\cap X\cong \A^1$ and the union 
$F_0\ \raisebox{1mm}{\rule{5mm}{0.2mm}}\ A$ is contractible to a smooth rational curve meeting the image of 
the component $M$.
\end{enumerate}  
Note that $X$ is an $\ML_0$-surface if and only if $A$ is a linear chain. We then call $X$ an {\em affine  pseudo-plane 
of $\ML_0$-type}. 
\svskip

If we are given a log deformation $(Y,\ol{Y},S,\ol{f},t_0)$ of the triple $(\ol{Y}_0,D_0,$ $Y_0)$, it follows 
by Ehresmann's fibration theorem that $p_g$ and the irregularity $q$ of the fiber $\ol{Y}_t$ is independent of $t$. 
Furthermore, by Lemma \ref{Lemma 2.2}, $Y_t$ has an $\A^1$-fibration if $Y_0$ has an $\A^1$-fibration. So, we can 
expect that $Y_t$ is an affine pseudo-plane if so is $Y_0$. Indeed, we have the following result. 

\begin{thm}\label{Theorem 4.6}
Let $\SF=(Y,\ol{Y},S,\ol{f},t_0)$ be a log deformation of $(\ol{Y}_0,D_0,$ $Y_0)$. Assume that $Y_0$ is an 
affine pseudo-plane. Then the following assertions hold.
\begin{enumerate}
\item[(1)]
$Y_t$ is an affine pseudo-plane for every point $t \in T$.
\item[(2)]
Assume that $Y_0$ is an affine pseudo-plane of $\ML_0$-type. Assume further that the boundary divisor $D_0$ in $\ol{Y}_0$ 
has the same weighted dual graph as above. Then $f : Y \to T$ is a trivial bundle with fiber $Y_0$ after shrinking $T$ 
if necessary.
\end{enumerate}
\end{thm}

\Proof
(1)~~We have only to show that $Y_t$ is an affine pseudo-plane for a small deformation of $Y_0$. After replacing 
$T$ by an \'etale finite covering, we may assume that $\ol{Y}_t$ is a minimal normal completion of $Y_t$ for every 
$t \in T$. Then, by Lemma \ref{Lemma 4.1}, the boundary divisor $D_t=S\cap \ol{Y}_t$ has the same weighted dual graph 
as shown above for $D_0$. Hence $Y_t$ has an $\A^1$-fibration of affine type. By Lemma \ref{Lemma 4.2}, 
$\Pic(Y_t) \cong \Pic(Y_0)$ which is a finite cyclic group. This implies that $Y_t$ is an affine pseudo-plane.

(2)~~Consider the completion $\ol{Y}_0$ of $Y_0$. We may assume that $\ol{Y}_0$ is a minimal normal completion of $Y_0$. 
In fact, a $(-1)$-curve contained in the boundary divisor $D_0$ which meets at most two other components of $D_0$ 
deforms to the nearby fibers and contracted simultaneously over the same $T$ by Lemma \ref{Lemma 2.5}, (2). 
Note that every fiber $Y_t$ has an $\A^1$-fibration of affine type by Lemma \ref{Lemma 2.2}. As in the proof of 
Lemma \ref{Lemma 2.2}, (2), by performing simultaneous (i.e., along the fibers of $\ol{f}$) blowing-ups and blowing-downs 
on the boundary $S$, we may assume that $Y_0$ has an $\A^1$-fibration which extends to a $\BP^1$-fibration on $\ol{Y}_0$ 
and that the boundary divisor $D_0$ has the weighted dual graph $\ell\ \raisebox{1mm}{\rule{5mm}{0.2mm}}\ M\ 
\raisebox{1mm}{\rule{5mm}{0.2mm}}\ A$ as specified in the condition (ii) above, where $A$ is a linear chain by the 
hypothesis. To perform a simultaneous blowing-up, we may have to choose as the center a cross-section on an irreducible 
component $S_i$ which is a $\BP^1$-bundle over $T$. If such a cross-section happens to intersect the curve $S_i\cap S_j$ 
with another component $S_j$, we shrink $T$ to avoid this intersection (see the remark in the second proof of Theorem 
\ref{Theorem 2.10}). Note that the interior $Y$ (more precisely, the inverse image of $f$ of the shrunken $T$) is not 
affected under these operations. Then the $(0)$ curve $\ell$ defines a $\BP^1$-fibration $\varphi : \ol{Y} \to V$ 
(see Lemma \ref{Lemma 2.5}, (1)). In particular, $\ell$ moves in an irreducible component, say $S_{-1}$, of $S$. 
The $(0)$ curve $M$ moves along the fibers of $\ol{f}$ in an irreducible component, say $S_0$, of $S$. 
By Lemma \ref{Lemma 4.1}, the curves in $A$ move along the fibers of $\ol{f}$ and fill out the irreducible components 
$S_1, \ldots, S_r$ of $S$. Hence $S=S_{-1}\cup S_0\cup S_1 \cup \cdots \cup S_r$ and $D_t=S\cdot\ol{Y}_t$ has the same 
weighted dual graph as $D_0$. 

Now consider a $(-1)$ curve $F_0$ on $\ol{Y}_0$. By Lemma \ref{Lemma 2.1}, $F_0$ moves along the fibers of $\ol{f}$ and 
fills out a smooth irreducible divisor $F$ which meets transversally an irreducible component $S_i\ (1 \le i \le r)$. 
In fact, the feather $F_0$ is unique on $Y_0$ and $(S_i\cdot F\cdot\ol{Y}_t)=(S_i\cdot F\cdot \ol{Y}_0)=1$. Let $S_1$ be 
the component of $S$ meeting $S_0$. Let $F_t=F\cap\ol{Y}_t$ and $S_{j,t}=S_j\cap\ol{Y}_t$ for every $t \in T$. Then 
$F_t+\sum_{j=2}^rS_{j,t}$ is contractible to a smooth point $P_t$ lying on $S_{1,t}$. After performing simultaneous 
elementary transformations on the fiber $\ell$ which is the fiber at infinity of the $\A^1$-fibration of the affine 
pseudo-plane $Y_t$, we may assume that $P_t$ is the intersection point $S_{0,t}\cap S_{1,t}$. By applying Lemma 
\ref{Lemma 2.5}, (2) repeatedly, we can contract $F$ and the components $S_2, \ldots, S_r$ simultaneously. Let $\ol{Z}$ 
be the threefold obtained from $\ol{Y}$ by these contractions. Then $\ol{Z}$ has a $\BP^1$-fibration 
$\psi : \ol{Z} \to V$ and the image of $S_0$ is a cross-section. Let $g=\sigma\cdot\psi : \ol{Z}\st{\psi} V \st{\sigma} T$ 
(see Lemma \ref{Lemma 2.1}, (3) for the notations). For every $t \in T$, $\ol{Z}_t:=g^{-1}(t)$ is a minimal $\BP^1$-bundle 
with a cross-section $S_{0,t}$. Since $(S_{0,t})^2=0$, $\ol{Z}_t$ is isomorphic to $\BP^1\times \BP^1$. Then $\ol{Z}$ is 
a trivial $\BP^1\times\BP^1$-bundle over $T$ after shrinking $T$ if necessary. In fact, $\ol{Z}$ with the images of $S_0$ 
and $S_{-1}$ removed is a deformation of $\A^2$, which is locally trivial in the Zariski topology by Theorem \ref{Theorem 2.10}. 
We may assume that $\psi : \ol{Z} \to V$ is the projection of $\BP^1\times\BP^1\times T$ onto the second and the third factors. 
Choose a section $\ol{S}'_0$ of $\psi$ which is disjoint from the image $\ol{S}_0$ of $S_0$. Then there is a non-trivial 
$G_m$-action on $\ol{Z}$ along the fibers of $\psi$ which has $\ol{S}_0$ and $\ol{S}'_0$ as the fixed point locus. 

Now reverse the contractions $\ol{Y} \to \ol{Z}$. The center of the first simultaneous blowing-up with center $S_0\cap S_1$
and the centers of the consecutive simultaneous blowing-ups except for the blowing-up which produces the component $F$ are 
$G_m$-fixed because the blowing-ups are fiberwise sub-divisional. Only the center $Q_t$ of the last blowing-up on $\ol{Y}_t$ 
is non-subdivisional. Let 
\[
\varphi : \ol{Y} \st{\sigma} \ol{Y}_1 \st{\sigma_1} \ol{Z}
\]
be the factorization of $\varphi$ where $\sigma$ is the last non-subdivisional blowing-up. By the construction, the natural 
$T$-morphism $\ol{f}_1 : \ol{Y}_1 \to T$ is a trivial fibration with fiber $(\ol{Y}_1)_0=\ol{f}_1^{-1}(t_0)$. Then there 
exists an element $\{\rho_t\}_{t \in T}$ of $G_m(T)$ such that $\rho_t(Q_{t_0})=Q_t$ for every $t \in T$ after shrinking 
$T$ if necessary. Here note that the $G_m$-action is nontrivial on the component with the point $Q_t$ thereon, for otherwise 
the $G_m$-action is trivial from the beginning. Then these $\{\rho_t\}_{t \in T}$ extends to a $T$-isomorphism 
$\wt{\rho} : \ol{Y}_0\times T \to \ol{Y}$, which induces a $T$-isomorphism $Y_0\times T \to Y$. Hence $Y$ is trivial.  
\QED

\section{Deformations of $\A^1$-fibrations of complete type}

In the setting of Theorem \ref{Theorem 2.8}, if the $\A^1$-fibration of a general fiber $Y_t$ is of complete type, we do 
not have the same conclusion. This case is treated in a recent work of Dubouloz and Kishimoto \cite{DK}. We consider this 
case by taking the same example of cubic surfaces in $\BP^3$ and explain how it is affine-uniruled.
\svskip

Taking a cubic hypersurface as an example, we first observe the behavior of the log Kodaira dimension for a flat family 
of smooth affine surfaces. Let $\stackrel{\vee}{\BP^3}$ be the dual projective $3$-space whose points correspond to 
the hyperplanes of $\BP^3$. We denote it by $T$. Let $S$ be a smooth cubic hypersurface in $\BP^3$ and let $\SW=S\times T$ 
which is a codimension one subvariety of $\BP^3\times T$. Let $\SH$ be the universal hyperplane in $\BP^3\times T$, 
which is defined by $\xi_0X_0+\xi_1X_1+\xi_2X_2+\xi_3X_3=0$, where $(X_0,X_1,X_2,X_3)$ and $(\xi_0,\xi_1,\xi_2,\xi_3)$ 
are respectively the homogeneous coordinates of $\BP^3$ and $T$. Let $\SD$ be the intersection of $\SW$ and $\SH$ in 
$\BP^3\times T$. Let $\pi : \SW \to T$ be the projection and let $\pi_\SD : \SD \to T$ be the restriction of $\pi$ onto $\SD$. 
Then $\pi$ and $\pi_\SD$ are the flat morphism. For a closed point $t \in T$, $\SW_t=\pi^{-1}(t)$ is identified with $S$ and 
$\SD_t=\pi_\SD^{-1}(t)$ is the hyperplane section $S\cap \SH_t$ in $\BP^3$, where $\SH_t$ is the hyperplane 
$\tau_0X_0+\tau_1X_1+\tau_2X_2+\tau_3X_3=0$ with $t=(\tau_0,\tau_1,\tau_2,\tau_3)$. Let $\SX=\SW\setminus\SD$ and 
$p : \SX \to T$ be the restriction of $\pi$ onto $\SX$. Then $\SX_t=p^{-1}(t)$ is an affine surface $S\setminus(S\cap \SH_t)$. 

Since $S$ is smooth, the following types of $S\cap \SH_t$ are possible. In the following, $F=0$ denotes the defining equation 
of $S$ and $H=0$ does the equation for $\SH_t$. 
\begin{enumerate}
\item[(1)]
A smooth irreducible plane curve of degree $3$.
\item[(2)]
An irreducible nodal curve, e.g., $F=X_0(X_1^2-X_2^2)-X_2^3+X_0^2X_3+X_3^3$ and $H=X_3$.
\item[(3)]
An irreducible cuspidal curve, e.g., $F=X_0X^2_1-X_2^3+X_3(X_0^2+X_1^2+X_2^2+X_3^2)$ amd $H=X_3$.
\item[(4)]
An irreducible conic and a line which meets in two points trans\-versally or in one point with multiplicity two. In fact,
let $\ell$ and $D$ be respectively a line and an irreducible conic in $\BP^2$ meeting in two points $Q_1, Q_2$, where 
$Q_1$ is possibly equal to $Q_2$. Let $C$ be a smooth cubic meeting $\ell$ in three points $P_i~(1 \le i \le 3)$ and $D$ 
in six points $P_i~(4 \le i \le 9)$, where the points $P_i$ are all distinct and different from $Q_1, Q_2$. Choose 
two points $P_1, P_2$ on $\ell$ and four points $P_i~(4 \le i \le 7)$ on $D$. Let $\sigma : S \to \BP^2$ be the blowing-up of 
these six points. Let $\ell', D'$ and $C'$ be the proper transforms of $\ell, D$ and $C'$. Then $S$ is a cubic 
hypersurface in $\BP^3$ and $K_S \sim -C'$. Since $\ell'+D' \sim C'$, it is a hyperplane section of $S$ with respect to 
the embedding $\Phi_{|C'|} : S \hookrightarrow \BP^3$. 
\item[(5)]
Three lines which are either meeting in one point or not. Let $\ell_i~(1 \le i \le 3)$ be the lines. Let $Q_1=\ell_1\cap \ell_3$ 
and $Q_2=\ell_2\cap \ell_3$. In the setting of (4) above, we consider $\ell=\ell_3$ and $D=\ell_1+\ell_2$. So, if $Q_1=Q_2$, 
three lines meet in one point. Choose a smooth cubic $C$ meeting three lines in nine distinct points $P_i~(1 \le i \le 9)$ 
other than $Q_1, Q_2$. Choose six points from the $P_i$, two points lying on each line. Then consider the blowing-up in 
these six points. The rest of the construction is the same as above. 
\end{enumerate}
Note that if $S$ is smooth $S\cap \SH_t$ cannot have a non-reduced component. In fact, the non-reduced component is a line 
in $\SH_t$. Hence we may write the defining equation of $S$ as 
\[
F=X_0^2(aX_1+X_0)+X_3G(X_0,X_1,X_2,X_3)=0,
\]
where $G=G(X_0,X_1,X_2,X_3)$ is a quadratic homogeneous polynomial and $a \in k$. We understand that $a=0$ if the non-reduced 
component has multiplicity three. By the Jacobian criterion, it follows that $S$ has singularities at the points $G=X_0=X_3=0$.

The affine surface $\SX_t$ has log Kodaira dimension $0$ in the cases (1), (2), (4) with the conic and the line meeting in 
two distinct points and (5) with non-confluent three lines, and $-\infty$ in the rest of the cases. Although $p : \SH \to T$ is 
a flat family of affine surfaces, the log Kodaira dimension drops to $-\infty$ exactly at the points $t \in T$ where the boundary 
divisor $S\cap \SH_t$ is not a divisor with normal crossings. This accords with a result of Kawamata concerning the invariance 
of log Kodaira dimension under deformations (cf. \cite{Kawamata}).

If $\lkd(\SX_t)=-\infty$, then $\SX_t$ has an $\A^1$-fibration. We note that if $\lkd(\SX_t)=0$ then $\SX_t$ has an 
$\A^1_*$-fibration. In fact, we consider the case where the boundary divisor $\SD_t$ is a smooth cubic curve. Then $S$ is 
obtained from $\BP^2$ by blowing up six points $P_i~(1 \le i \le 6)$ on a smooth cubic curve $C$. Choose four points 
$P_1,P_2,P_3,P_4$ and let $\Lambda$ be a linear pencil of conics passing through these four points. Let $\sigma : S \to \BP^2$ 
be the blowing-up of six points $P_i~(1 \le i \le 6)$. The proper transform $\sigma'\Lambda$ defines a $\BP^1$-fibration 
$f : S \to \BP^1$ for which the proper transform $C'=\sigma'(C)$ is a $2$-section. Since $\SX_t$ is isomorphic to $S\setminus C'$, 
$\SX_t$ has an $\A^1_*$-fibration. 

Looking for an $\A^1$-fibration in the case $\lkd(\SX_t)=-\infty$ is not an easy task. Consider, for example, the case where 
$X=\SX_t$ is obtained as $S\setminus(Q\cup \ell)$, where $Q$ is a smooth conic and $\ell$ is a line in $\BP^2$ which meet in 
one point with multiplicity two. As explained in the above, such an $X$ is obtained from $\BP^2$ by blowing up six points 
$P_1, \ldots, P_6$ such that $P_1, P_2$ lie on a line $\wt{\ell}$ and $P_3,P_4,P_5,P_6$ are points on a conic $\wt{Q}$. Then 
the proper transforms on $S$ of $\wt{\ell}, \wt{Q}$ are $\ell, Q$. Consider the linear pencil $\wt{\Lambda}$ on $\BP^2$ spanned 
by $2\wt{\ell}$ and $\wt{Q}$. Then a general member of $\Lambda$ is a smooth conic meeting $\wt{Q}$ in one point 
$\wt{Q}\cap\wt{\ell}$ with multiplicity four. The proper transform $\Lambda$ of $\wt{\Lambda}$ on $S$ defines an $\A^1$-fibration 
on $X$.
\svskip

The following result of Dubouloz-Kishimoto except for the assertion (4) was orally communicated to one of the authors (see 
\cite{DK}).

\begin{thm}\label{Theorem 5.1}
Let $S$ be a cubic hypersurface in $\BP^3$ with a hyperplane section $S\cap H$ which consists of a line and a conic meeting 
in one point with multiplicity two. Let $Y=\BP^3\setminus S$ which is a smooth affine threefold. Then the following assertions hold.
\begin{enumerate}
\item[(1)]
$\lkd(Y)=-\infty$.
\item[(2)]
Let $f : Y \to \A^1$ be a fibration induced by the linear pencil on $\BP^3$ spanned by $S$ and $3H$. Then a general fiber 
$Y_t$ of $f$ is a cubic hypersurface $S_t$ minus $Q\cup \ell$, where $Q$ is a conic and $\ell$ is a line which meet 
in one point with multiplicity two. Hence $\lkd(Y_t)=-\infty$ and $Y_t$ has an $\A^1$-fibration. 
\item[(3)]
$Y$ has no $\A^1$-fibration.
\item[(4)]
There is a finite covering $T'$ of $\A^1$ such that the normalization of $Y\times_{\A^1}T'$ has an $\A^1$-fibration.
\end{enumerate}
\end{thm}
\Proof
(1)\ Since $K_{\BP^3}+S \sim -4H+3H=-H$, it follows that $\lkd(Y)=-\infty$.

(2)\ The pencil spanned by $S$ and $3H$ has base locus $Q\cup\ell$ and its general member, say $S_t$, is a cubic hypersurface 
containing $Q\cup\ell$ as a hyperplane section. It is clear that $S_t\setminus(Q\cup\ell)=Y_t$. Hence, as explained above, 
$Y_t$ has an $\A^1$-fibration.

(3)\ Let $\tau : \wt{S} \to \BP^3$ be the cyclic triple covering of $\BP^3$ ramified totally over the cubic hypersurface $S$. 
Then $\wt{S}$ is a cubic hypersurface in $\BP^4$ and $\tau^*(S)=3\wt{H}$, where $\wt{H}$ is a hyperplane in $\BP^4$. The 
restriction of $\tau$ onto $Z:=\wt{S}\setminus \wt{S}\cap\wt{H}$ induces a finite \'etale covering $\tau_Z : Z \to Y$. 
Suppose that $Y$ has an $\A^1$-fibration $\varphi : Y \to T$. Then $T$ is a rational surface. Since $\tau_Z$ is finite \'etale, 
this $\A^1$-fibration $\varphi$ lifts up to an $\A^1$-fibration $\wt{\varphi} : Z \to \wt{T}$. By \cite{CG}, $\wt{S}$ is 
unirational and irrational. Hence $\wt{T}$ is a rational surface. This implies that $Z$ is a rational threefold. This 
is a contradiction because $\wt{S}$ is irrational.

(4)\ There is an open set $T$ of $\A^1$ such that the restriction of $f$ onto $f^{-1}(T)$ is a smooth morphism onto $T$. 
By abuse of the notations, we denote $f^{-1}(T)$ by $Y$ anew and the restriction of $f$ onto $f^{-1}(T)$ by $f$. Hence 
$f : Y \to T$ is a smooth morphism. Let $K=k(t)$ be the function field of $T$ and let $Y_K$ be the generic fiber. Let 
$\ol{K}$ be an algebraic closure of $K$. Then $Y_{\ol{K}}:=Y_K\otimes_K\ol{K}$ is identified with 
$S_{\ol{K}}\setminus (Q\cup\ell)$, where $S_{\ol{K}}$ is a cubic hypersurface in $\BP^3_{\ol{K}}$ defined by 
$F_K=F_0+tX_3^3=0$. Here $t$ is a coordinate of $\A^1$ and $(X_0,X_1,X_2,X_3)$ is a system of homogeneous coodinates of 
$\BP^3$ such that $F_0(X_0,X_1,X_2,X_3)=0$ is the defining equation of the cubic hypersurface $S$ and the hyperplane $H$ 
is defined by $X_3=0$. Then $Y_{\ol{K}}$ is obtained from $\BP^2_{\ol{K}}$ by blowing up six $\ol{K}$-rational points 
in general position (two points on the image of $\ell$ and four points on the image of $Q$). As explained earlier, 
there is an $\A^1$-fibration on $Y_{\ol{K}}$ which is obtained from conics on $\BP^2_{\ol{K}}$ belonging to the pencil 
spanned by $Q$ and $2\ell$. This construction involves six points on $\BP^2_{\ol{K}}$ to be blown up to obtain the cubic 
hypersurface $S_{\ol{K}}$ and four points (the point $Q\cap\ell$ and its three infinitely near points). Hence there 
exists a finite algebraic extension $K'/K$ such that all these ponts are rational over $K'$. Let $T'$ be the normalization 
of $T$ in $K'$. Let $Y'=Y\otimes_KK'$. Then $Y'$ has an $\A^1$-fibration. 
\QED

Based on the assertion (4) above, we propose the following conjecture.

\begin{conj}\label{Conjecture 5.2}
Let $f : Y \to T$ be a smooth morphism from a smooth affine threefold $Y$ onto a smooth affine curve $T$ such that every 
closed fiber $Y_t$ has an $\A^1$-fibration of complete type. Then there exists a finite covering $T'$ of $T$ such that 
the normalization of $Y\times_TT'$ has an $\A^1$-fibration.
\end{conj}

\begin{remark}\label{Remark 5.3}{\em
The conjecture \ref{Conjecture 5.2} is true if Theorem \ref{Theorem 2.8} holds after an \'etale finite base change of 
$T$ in the case where the general fibers of $f$ have $\A^1$-fibrations of complete type. A main obstacle in trying 
to extend the proof in the case of $\A^1$-fibrations of affine type is to show that, with the notations in the proof of 
Theorem \ref{Theorem 2.8}, the locus of base points $P_t$ of the linear pencil $\Lambda_t$ on $\ol{Y}_t$ with $t$ 
varying in an open neighborhood of $t_0 \in T$ (resp. the loci of infinitely near base points) is a cross-section of the 
$\BP^1$-bundle $\ol{f}|_{S_1} : S_1 \to T$ (resp. the exceptional $\BP^1$-bundle).}
\end{remark}


\begin{thebibliography}{30}

\bibitem{B}
G. Bredon, Topology and Geometry, Graduate Texts in Mathematics no. 139, Springer.

\bibitem{CG}
C.H. Clemens and P.A. Griffiths, The intermediate Jacobian of the cubic threefold. Ann. of Math. (2) {\bf 95} (1972), 281--356. 

\bibitem{DK}
A. Dubouloz and T. Kishimoto, Log-uniruled affine varieties without cylinder-like open subsets, arXiv: 1212.0521, 2012.

\bibitem{Fujita}
T. Fujita, On the topology of noncomplete algebraic surfaces, J. Fac. Sci. Univ. Tokyo Sect. IA Math. {\bf 29} (1982), no. 3, 
503--566.

\bibitem{FZ}
H. Flenner and M. Zaidenberg, $\Q$ -acyclic surfaces and their deformations, Classification of algebraic varieties (L'Aquila, 1992), 
143--208, Contemp. Math. {\bf 162}, Amer. Math. Soc., 1994.

\bibitem{FKZ}
H. Flenner, Sh. Kaliman and M. Zaidenberg, Deformation equivalence of affine ruled surfaces, arXiv:1305.5366v1, 23 May 2013.

\bibitem{G}
A. Grothendieck, Techniques de construction et th\'eor\`emes d'existence en g\'eom\'etrie alg\'ebrique IV: Les sch\'emas 
de Hilbert, S\'eminaire Bourbaki, 13e ann\'ee, 1960/61, no. 221, 1--28. 

\bibitem{GKM}
R.V. Gurjar, K. Masuda and M. Miyanishi, $\A^1$-fibrations on affine threefolds, J. pure and applied algebra 
{\bf 216} (2012), 296--313. 

\bibitem{GKM2}
R.V. Gurjar, K. Masuda and M. Miyanishi, Surjective derivations in small dimensions, preprint, to appear in the commemorative 
volume to the eightieth birthday of C.S. Seshadri.

\bibitem{GMMR}
R.V. Gurjar, K. Masuda, M. Miyanishi and P. Russell, Affine lines on affine surfaces and the Makar-Limanov invariant, 
Canad. J. Math. {\bf 60} (1), 2008, 109--139.

\bibitem{H}
R. Hartshorne, Algebraic geometry, Springer, 1977.

\bibitem{Iitaka}
S. Iitaka, Deformations of compact complex surfaces, Global Analysis (Papers in Honor of K. Kodaira), 267--272, 
Univ. Tokyo Press, Tokyo, 1969.  

\bibitem{Kamb}
T. Kambayashi, On the absence of nontrivial separable forms of the affine plane. J. Algebra {\bf 35} (1975), 
449--456. 

\bibitem{KM}
T. Kambayashi and M. Miyanishi, On flat fibrations by the affine line, Illinois J. Math. {\bf 22}
(1978), 662--671.

\bibitem{KW}
T. Kambayashi and D. Wright, Flat families of affine lines are affine line bundles, Illinois J. Math. {\bf 29} 
(1985), 672--681.

\bibitem{KZ}
Sh. Kaliman and M. Zaidenberg, Families of affine planes: the existence of a cylinder, Michigan Math. J. {\bf 49} (2001), 
no. 2, 353--367.

\bibitem{Kawamata1}
Y. Kawamata, On deformations of compactifiable complex manifolds, Math. Ann. {\bf 235} (1978), no. 3, 247--265.

\bibitem{Kawamata}
Y. Kawamata, On the classification of non-complete algebraic surfaces, Algebraic geometry 
(Proc. Summer Meeting, Univ. Copenhagen, Copenhagen, 1978), 215--232, Lecture Notes in Math., {\bf 732}, 
Springer, Berlin, 1979. 

\bibitem{K1}
K. Kodaira, On the stability of compact submanifolds of complex mnaifolds, Amer. J. Math. {\bf 85} (1963), 79--94.

\bibitem{K}
K. Kodaira, Complex manifolds and deformation of complex structures, Classics in Mathematics, Springer, 1986.

\bibitem{Kollar}
J. Kollar, Rational curves on algebraic varieties, Ergebnisse der Mathematik und ihrer Grenzgebiete {\bf 32}, 
Springer-Verlag, Berlin, 1996.

\bibitem{KollarMori}
J. Kollar and S. Mori, Birational geometry of algebraic varieties, Cambridge Tracts in Mathematics {\bf 134}, 
Cambridge University Press, Cambridge, 1998. 

\bibitem{M1}
M. Miyanishi, An algebro-topological characterization of the affine space of dimension three, Amer. J. Math. 
{\bf 106} (1984), 1469-1486.

\bibitem{M2}
M. Miyanishi, Affine pseudo-coverings of algebraic surfaces, J. Algebra {\bf 294} (2005), 156--176.

\bibitem{Morrow}
J.A. Morrow, Minimal normal compactifications of $\C^2$, Complex Analysis, 1972 (Proc. Conf. Rice Univ., Houston, Tex.,
1972), Vol. I: Geometry of singularities, Rice Univ. Studies {\bf 59} (1973), no. 1, 97--112.

\bibitem{N}
W.D. Neumann, On the topology of curves in complex surfaces, Topological methods in algebraic transformation groups 
(New Brunswick, NJ, 1988), 117--133, Progr. Math., {\bf 80}, Birkhauser Boston, Boston, MA, 1989.

\bibitem{R}
C.P. Ramanujam, A topological characterization of the affine plane as an algebraic variety, Ann. of Math. {\bf 94} 
(1971), 69--88.

\bibitem{Russell}
P. Russell, Embeddings problems in affine algebraic geometry, Polynomial automorphisms and related topics,
113--135, Publishing House for Science and Technology, Hanoi, 2007.

\bibitem{S}
A. Sathaye, Polynomial ring in two variables over a DVR: a criterion, Invent. Math. {\bf 74} (1983), no. 1, 159--168. 

\bibitem{W}
R.O. Wells, Differential Analysis on Complex Manifolds, Graduate Texts in Mathematics no. 65, 
Springer.

\end{thebibliography}
\end{document}